\documentclass{amsart}

\title{Reflective algebras over factorizable Hopf algebras}
\author[R. Sugitani]{Rei Sugitani}
\address[R. Sugitani]{Graduate School of Pure and Applied Sciences,
	University of Tsukuba,
	1-1-1 Tennodai, Tsukuba, Ibaraki 305-8571, Japan}
\email{rsugitani@math.tsukuba.ac.jp}

\usepackage{amsmath,amssymb,amsthm,mathtools}
\usepackage{bbm,dsfont,mathrsfs}
\usepackage{hyperref}
\usepackage{cleveref}
\usepackage{enumitem}

\numberwithin{equation}{section}

\newtheorem{theorem}{Theorem}[section]
\newtheorem{lemma}[theorem]{Lemma}
\newtheorem{corollary}[theorem]{Corollary}
\newtheorem{proposition}[theorem]{Proposition}

\theoremstyle{definition}
\newtheorem{definition}[theorem]{Definition}

\theoremstyle{remark}
\newtheorem{remark}[theorem]{Remark}
\newtheorem{example}[theorem]{Example}

\crefname{section}{Section}{Sections}
\crefname{subsection}{Section}{Sections}

\crefname{lemma}{Lemma}{Lemmas}
\crefname{theorem}{Theorem}{Theorems}
\crefname{corollary}{Corollary}{Corollaries}
\crefname{proposition}{Proposition}{Propositions}
\crefname{definition}{Definition}{Definitions}
\crefname{remark}{Remark}{Remarks}
\crefname{example}{Example}{Examples}
\crefname{hypothesis}{Hypothesis}{Hypotheses}

\newcommand{\bfk}{\Bbbk}
\newcommand{\id}{\mathrm{id}}
\newcommand{\op}{\mathrm{op}}
\newcommand{\cop}{\mathrm{cop}}
\newcommand{\Hom}{\mathrm{Hom}}
\newcommand{\co}{\mathrm{co}}

\newcommand{\HAlg}{{}^H\!\mathrm{Alg}}

\DeclareMathOperator{\Ker}{Ker}
\DeclareMathOperator{\Img}{Im}

\begin{document}

\begin{abstract}
	The reflective algebra $(R_H(A),K_{\mathrm{ref}}(A))$, introduced by Laugwitz, Walton and Yakimov \cite{LWY25}, is a quasitriangular left $H$-comodule algebra constructed from a quasitriangular Hopf algebra $(H,R)$ and a left $H$-comodule algebra $A$.
	We study the structure and factorizability of finite-dimensional reflective algebras.
	When $(H,R)$ is factorizable, we show that $(R_H(A),K_{\mathrm{ref}}(A))$ is isomorphic, as a quasitriangular left $H$-comodule algebra, to $A\otimes H$ equipped with the trivial left $H$-coaction on $A$ and the universal $K$-matrix induced by $R_{21}R$.
	In this case, we also classify all universal $K$-matrices on $R_H(A)$.
	Finally, we prove that $(R_H(A),K_{\mathrm{ref}}(A))$ is factorizable if and only if $(H,R)$ is factorizable and $A$ is simple as an algebra.
\end{abstract}

\maketitle

\section{Introduction}
\label{sec:introduction}

Braided tensor categories play an important role in various areas of mathematics and mathematical physics, including representation theory, low-dimensional topology, and quantum field theory \cite{EGNO15,RT91,FRS02}.
Brochier \cite{Bro13} introduced braided module categories, which are module categories over braided tensor categories equipped with an additional braiding.
Braided module categories provide a categorical framework for studying invariants of tangles in a solid torus \cite{Bro13}.
They also arise in the theory of quantum character varieties and extensions of braided tensor categories \cite{BZBJ18,DN21}.

On the Hopf algebra side, quasitriangular comodule algebras provide algebraic realizations of braided module categories.
Let $(H,R)$ be a quasitriangular Hopf algebra and $B$ a left $H$-comodule algebra.
A quasitriangular left $H$-comodule algebra is a pair $(B,K)$, where $K\in H\otimes B$ is a universal $K$-matrix \cite{Kol20}.
If $B$ is quasitriangular, the category of left $B$-modules is then a braided module category over the braided monoidal category of left $H$-modules.
Universal $K$-matrices provide algebraic descriptions of module braidings and satisfy the quantum reflection equation.
They also give rise to representations of braid groups of type $B$ and appear in integrable models with boundary conditions \cite{Kol20,Skl88}.
In view of these applications, the construction of quasitriangular comodule algebras and the explicit description and classification of universal $K$-matrices are natural problems.
Indeed, $K$-matrices associated with quantum symmetric pairs have been studied extensively \cite{BK19,AV22,AV25}.
For finite-dimensional quasitriangular Hopf algebras, $K$-matrices on certain coideal subalgebras have been classified \cite{MW26}.

To provide general constructions of braided module categories and quasitriangular comodule algebras, Laugwitz, Walton, and Yakimov \cite{LWY25} introduced reflective centers and reflective algebras.
For a braided tensor category $\mathcal C$ and a $\mathcal C$-module category $\mathcal M$, the reflective center $\mathcal E_{\mathcal C}(\mathcal M)$ is a braided $\mathcal C$-module category.
Given a finite-dimensional quasitriangular Hopf algebra $(H,R)$ and a finite-dimensional left $H$-comodule algebra $A$, one can construct the reflective algebra $R_H(A)$, which is a quasitriangular left $H$-comodule algebra equipped with the canonical universal $K$-matrix $K_{\mathrm{ref}}(A)$.
The category $R_H(A)\text{-}\mathrm{mod}$ of finite-dimensional left $R_H(A)$-modules is equivalent to $\mathcal E_{H\text{-}\mathrm{mod}}(A\text{-}\mathrm{mod})$ as a braided module category over $H\text{-}\mathrm{mod}$, where $H\text{-}\mathrm{mod}$ and $A\text{-}\mathrm{mod}$ denote the categories of finite-dimensional left $H$-modules and left $A$-modules, respectively.
This construction parallels the realization of the Drinfeld center by the Drinfeld double.
In this way, quasitriangular left $H$-comodule algebras can be constructed systematically.

Besides the construction of quasitriangular comodule algebras, a natural question is when the corresponding braided module categories are non-degenerate.
For braided finite tensor categories, non-degeneracy is known to be equivalent to factorizability \cite{Shi19}.
Walton and Yadav \cite{WY25} extended this result to braided module categories by defining non-degeneracy and factorizability for exact and indecomposable braided module categories and proving their equivalence.
For a finite-dimensional quasitriangular left $H$-comodule algebra $(B,K)$, they also introduced the notion of factorizability under the assumption that $B\text{-}\mathrm{mod}$ is an exact and indecomposable module category over $H\text{-}\mathrm{mod}$.
They then proved that $(B,K)$ is factorizable if and only if $B\text{-}\mathrm{mod}$ is non-degenerate.
For reflective algebras, they further showed that $(R_H(A),K_{\mathrm{ref}}(A))$ is factorizable whenever $R_H(A)$ is $H$-simple \cite{WY25}.
Determining when reflective algebras are $H$-simple is therefore important for studying the non-degeneracy of reflective centers.

The reflective algebra $R_H(A)$ is constructed in \cite{LWY25} as a crossed product of $A$ with the dual of an $H$-module coalgebra determined by $R$, with explicit formulas for its multiplication and left $H$-coaction.
However, this construction does not always make it easy to describe the structure or $H$-simplicity of $R_H(A)$ directly in terms of properties of $H$ and $A$.
This raises the question of which conditions on $H$ and $A$ characterize the $H$-simplicity of $R_H(A)$ and the factorizability of $(R_H(A),K_{\mathrm{ref}}(A))$.
The aim of this paper is to address this question and to study the structure and universal $K$-matrices of reflective algebras.

Our study is based on the fact that the Drinfeld map $Q_H:H^*\to H$ can be regarded as a morphism of quasitriangular left $H$-comodule algebras.
More precisely, identifying $R_H(\bfk)$ with $H^*$ as vector spaces and equipping $H$ with its coproduct as a left $H$-coaction and with the universal $K$-matrix $R_{21}R$, we obtain the morphism $Q_H:(R_H(\bfk),K_{\mathrm{ref}}(\bfk))\to(H,R_{21}R)$ (\cref{prop:drinfeld-map-QT-map}).
We first realize the quotient comodule algebra $R_H(A)/A\,\Ker Q_H$ as a cotensor product of $H$ and $A$ over a quotient coalgebra associated with $\Img Q_H$ (\cref{thm:cotensor-quasitriangular}).
We also show that this quotient is $H$-simple if and only if $A$ is simple with respect to the induced coaction of this quotient coalgebra (\cref{cor:quotient-H-simple}).
Next, suppose that $(H,R)$ is factorizable, or equivalently that $Q_H$ is an isomorphism of vector spaces.
We describe the structure of a general quasitriangular left $H$-comodule algebra $(B,K)$ in terms of its coinvariant subalgebra $B^{\co H}$ and investigate its universal $K$-matrices, $H$-simplicity, and factorizability (\cref{sec:factorizable-case}).
Applying these results to reflective algebras, we obtain an isomorphism of quasitriangular left $H$-comodule algebras
\begin{equation}\label{eq:main-result-1}
(R_H(A),K_{\mathrm{ref}}(A))\cong(A_{\mathrm{triv}}\otimes H,(R_{21}R)_{13}),
\end{equation}
as stated in \cref{thm:factorizable-reflective-structure}, where $A_{\mathrm{triv}}$ denotes the algebra $A$ equipped with the trivial left $H$-coaction.
Using this isomorphism, we show that the universal $K$-matrices on $R_H(A)$ are in bijection with algebra maps from $H$ to the center $Z(A)$ of $A$.
We also show that the $H$-costable ideals of $R_H(A)$ correspond to the ideals of $A$, so that $R_H(A)$ is $H$-simple if and only if $A$ is simple as an algebra.
Combining these results, we prove that for a finite-dimensional quasitriangular Hopf algebra $(H,R)$ and a nonzero finite-dimensional left $H$-comodule algebra $A$, the reflective algebra $(R_H(A),K_{\mathrm{ref}}(A))$ is factorizable if and only if $(H,R)$ is factorizable and $A$ is simple as an algebra (\cref{thm:reflective-factorizability-criterion}).
Finally, we take $H$ to be the small quantum group $u_q(\mathfrak{sl}_2)$ and consider the reflective algebras associated with $H$ itself and with the left $H$-comodule algebra $A_2(r;\zeta)$ introduced in \cite{NSS25}.
We give explicit descriptions of the corresponding reflective algebras and their universal $K$-matrices, and determine their $H$-simplicity and factorizability.

The isomorphism \eqref{eq:main-result-1} may be viewed as an analogue for reflective algebras of the well-known result that the Drinfeld double $D(H)$ of a factorizable Hopf algebra $H$ is isomorphic to $H\otimes H$ as an algebra (see, for example, \cite[Theorem~13.6.1]{Rad12}).
Moreover, if $\bfk$ is algebraically closed and $A$ is simple as an algebra, \eqref{eq:main-result-1} implies that $R_H(A)\text{-}\mathrm{mod}$ is equivalent, as a braided module category over $H\text{-}\mathrm{mod}$, to the regular braided module category $H\text{-}\mathrm{mod}$.
This is analogous to the result that every indecomposable braided module category over a non-degenerate braided fusion category is equivalent to the regular braided module category \cite[Proposition~4.17]{DN21}.
Thus the reflective center $\mathcal E_{H\text{-}\mathrm{mod}}(A\text{-}\mathrm{mod})$ satisfies a similar equivalence without the assumption that $H$ is semisimple.

\subsection*{Organization of this paper}

This paper is organized as follows:
In \cref{sec:preliminaries}, we fix notation and recall basic facts about Hopf algebras and comodule algebras.
We first recall the definitions and basic properties of quasitriangular Hopf algebras and quasitriangular left $H$-comodule algebras.
We then introduce cotensor products and relative Hopf modules and present Takeuchi's equivalence \cite{Tak79}.
We also discuss the relationship between the $H$-simplicity of a finite-dimensional left $H$-comodule algebra and the exactness and indecomposability of its module category of finite-dimensional modules.

In \cref{sec:reflective-centers-algebras}, we recall known results on reflective centers and reflective algebras and establish the properties needed in this paper.
We first review the multiplication, left $H$-coaction, and canonical universal $K$-matrix of the reflective algebra $(R_H(A),K_{\mathrm{ref}}(A))$.
Using the universal property of $R_H(\bfk)$, we then consider the morphism of quasitriangular left $H$-comodule algebras from $(R_H(\bfk),K_{\mathrm{ref}}(\bfk))$ to $(H,R_{21}R)$ and show that it coincides with the Drinfeld map $Q_H$ under the identification of $R_H(\bfk)$ with $H^*$ as vector spaces.
Finally, we recall the definition of factorizability for quasitriangular comodule algebras introduced by Walton and Yadav, together with related results.

In \cref{sec:drinfeld-quotient}, we study the kernel and image of the Drinfeld map $Q_H$ and use them to construct a quotient of the reflective algebra.
We first show that $\Ker Q_H$ is stable under the right $H$-action.
This implies that the subspace $A\,\Ker Q_H\subset R_H(A)$ spanned by products of elements of $A$ and $\Ker Q_H$ is an $H$-costable ideal of $R_H(A)$, and we equip the quotient with a universal $K$-matrix.
We then define a quotient Hopf algebra $T$ associated with $\Img Q_H$ and use Takeuchi's equivalence to realize the quotient as the cotensor product $H\square_T A$.
We also characterize the $H$-simplicity of the quotient in terms of the $T$-simplicity of $A$.

In \cref{sec:factorizable-case}, we assume that $(H,R)$ is factorizable and study a general quasitriangular left $H$-comodule algebra $(B,K)$.
The Drinfeld map allows us to define a right $H$-action on $B$ that makes $B$ a Hopf module.
We apply the fundamental theorem of Hopf modules to obtain an isomorphism of quasitriangular left $H$-comodule algebras $(B^{\co H}\otimes H,(R_{21}R)_{13})\cong(B,K)$.
We next describe the universal $K$-matrices on $B$ using the relationship between algebra maps from $H$ to $B^{\co H}$ and maps of left $H$-comodule algebras.
Finally, we characterize the bijectivity of the map $\theta_B$ appearing in the definition of a factorizable comodule algebra \cite{WY25} and clarify its relation to the factorizability of $(B,K)$.

In \cref{sec:factorizable-reflective-case}, we apply the general results of \cref{sec:factorizable-case} to reflective algebras using the results of \cref{sec:drinfeld-quotient}.
When $(H,R)$ is factorizable, the description of the quotient obtained earlier gives an explicit isomorphism of quasitriangular left $H$-comodule algebras for $R_H(A)$ and identifies its coinvariant subalgebra with $A$.
The results of the preceding section then yield a classification of universal $K$-matrices on $R_H(A)$ and show that $\theta_{R_H(A)}$ is bijective if and only if $Z(A)=\bfk$.
We finally remove the factorizability assumption on $H$ and use the $H$-costable ideal constructed in \cref{sec:drinfeld-quotient} to characterize the $H$-simplicity of $R_H(A)$ and derive necessary and sufficient conditions for the factorizability of $(R_H(A),K_{\mathrm{ref}}(A))$.

In \cref{sec:example-uqsl2}, we take $H$ to be the small quantum group $u_q(\mathfrak{sl}_2)$ associated with a root of unity $q$ of odd order and give examples of the above results.
We first recall the quasitriangular structure and factorizability of $H$ and describe $R_H(\bfk)$.
For the left $H$-comodule algebra $A_2(r;\zeta)$ introduced in \cite{NSS25}, we then explicitly determine the algebra structure, left $H$-coaction, and universal $K$-matrix of its reflective algebra.
We also determine $H$-simplicity and factorizability in these examples.

\subsection*{Acknowledgments}
The author thanks Kenichi Shimizu for helpful comments.
The author thanks Taiki Shibata for pointing out an error in Section~2.3 of an earlier version of this paper.
The author thanks Akira Masuoka for valuable comments and suggestions, which helped improve the clarity and presentation of this paper.

\section{Preliminaries}
\label{sec:preliminaries}

\subsection{Conventions}
\label{subsec:hopf-comodule}

Throughout this paper, we work over an algebraically closed field $\bfk$.
In \cref{sec:example-uqsl2}, we further assume that $\bfk$ has characteristic zero.
For a $\bfk$-vector space $V$, we denote its dual space by $V^*$ and the canonical pairing of $f\in V^*$ and $v\in V$ by $\langle f,v\rangle$.
The symbol $\otimes$ denotes the tensor product over $\bfk$.
By an algebra or a coalgebra, we mean an associative unital algebra or a coassociative counital coalgebra over $\bfk$, respectively.
Unless otherwise stated, all algebras, coalgebras, Hopf algebras, modules, and comodules are finite-dimensional.
The comultiplication and counit of a coalgebra $C$ are denoted by $\Delta:C\to C\otimes C$ and $\varepsilon:C\to\bfk$, respectively.
We use Sweedler's notation for comultiplication:
\begin{equation*}
\Delta(c)=c_{(1)}\otimes c_{(2)},\quad
\Delta(c_{(1)}) \otimes c_{(2)} = c_{(1)}\otimes c_{(2)}\otimes c_{(3)} = c_{(1)} \otimes \Delta (c_{(2)}).
\end{equation*}
Let $C$ be a coalgebra.
For a left $C$-comodule $M$, we denote its coaction by $\delta_M:M\to C\otimes M$ and write $\delta_M(m)=m_{(-1)}\otimes m_{(0)}$ for $m\in M$, following Sweedler's notation.
Similarly, the coaction of a right $C$-comodule $N$ is denoted by $\rho_N:N\to N\otimes C$, with $\rho_N(n)=n_{(0)}\otimes n_{(1)}$.
When we simply write $C$-comodule without specifying left or right, we mean a left $C$-comodule.

We refer the reader to the books by Montgomery and Kassel \cite{Mon93,Kas95} for the basics of Hopf algebras.
For a Hopf algebra $H$, we denote its antipode by $S$.
The antipode of a finite-dimensional Hopf algebra is bijective.
We denote by $H^{\cop}$ the Hopf algebra with the same algebra structure as $H$, comultiplication $\Delta^{\op}:=\tau_{H,H}\circ\Delta$, and antipode $S^{-1}$, where $\tau_{H,H}:H\otimes H\to H\otimes H$ is the flip map defined by $x\otimes y\mapsto y\otimes x$.
For a left $H$-comodule algebra $B$, we denote its coinvariant subalgebra by $B^{\co H}=\{b\in B\mid\delta_B(b)=1_H\otimes b\}$.

We denote the category of finite-dimensional left $H$-modules by $H\text{-}\mathrm{mod}$.
We also denote the category of finite-dimensional left $H$-comodules by ${}^{H}\mathrm{comod}$.
Both are finite tensor categories.

\subsection{Quasitriangular comodule algebras}
\label{subsec:quasitriangular}

An \emph{universal quantum $R$-matrix} for a Hopf algebra $H$ is an invertible element $R\in H\otimes H$ satisfying
\begin{equation*}
(\Delta\otimes\id_H)(R)=R_{13}R_{23},\quad (\id_H\otimes\Delta)(R)=R_{13}R_{12},\quad R\Delta(h)=\Delta^{\op}(h)R
\end{equation*}
for all $h\in H$.
Writing $R=\sum_i s_i\otimes t_i$, we use $R_{ab}$ to denote the element of $H^{\otimes 3}$ obtained by placing $s_i$ and $t_i$ in the $a$-th and $b$-th tensor factors, respectively, and $1_H$ in the remaining tensor factor.
For example, $R_{13}=\sum_i s_i\otimes 1_H\otimes t_i$.
We also write $R_{21}=\sum_i t_i\otimes s_i$.
A Hopf algebra $H$ equipped with a quantum $R$-matrix $R$ is called a \emph{quasitriangular Hopf algebra} and is denoted by $(H,R)$.
Throughout this subsection, let $(H,R)$ be a quasitriangular Hopf algebra.

\begin{definition}
\label{def:drinfeld-map-factorizable}
The \emph{Drinfeld map} is the linear map $Q_H:H^*\to H$ defined by
\begin{equation*}
Q_H(f)=(f\otimes\id_H)(R_{21}R)
\end{equation*}
for $f\in H^*$.
The quasitriangular Hopf algebra $(H,R)$ is called \emph{factorizable} if $Q_H:H^*\to H$ is a linear isomorphism.
\end{definition}

\begin{definition}
	\label{def:quantum-K-matrix}
	An \emph{universal $K$-matrix} on an $H$-comodule algebra $B$ is an invertible element $K\in H\otimes B$ satisfying
	\begin{align}
		\label{eq:K-coproduct}
		(\Delta\otimes\id_B)(K)
		&=
		K_{23}R_{21}K_{13}R_{21}^{-1}
		\\
		\label{eq:K-coaction}
		(\id_H\otimes\delta_B)(K)
		&=
		R_{21}K_{13}R_{12}
		\\
		\label{eq:K-commutative}
		K\delta_B(b)
		&=
		\delta_B(b)K
		\quad
		(b\in B).
	\end{align}
	Here, if $K=\sum_i g_i\otimes p_i$, then
	\begin{equation*}
		K_{13}=\sum_i g_i\otimes1_H\otimes p_i \in H\otimes H\otimes B,
		\quad
		K_{23}=\sum_i1_H\otimes g_i\otimes p_i \in H\otimes H\otimes B.
	\end{equation*}
	A left $H$-comodule algebra $B$ equipped with a universal $K$-matrix $K$ is called a \emph{quasitriangular left $H$-comodule algebra} and is denoted by $(B,K)$.
\end{definition}

\begin{proposition}
\label{prop:K-matrix-quotient}
Let $(B,K)$ be a quasitriangular left $H$-comodule algebra and $C$ a left $H$-comodule algebra.
For a surjective left $H$-comodule algebra map $f:B\to C$, the element
\begin{equation*}
K_C:=(\id_H\otimes f)(K)\in H\otimes C
\end{equation*}
is a universal $K$-matrix on $C$.
In particular, $(C,K_C)$ is a quasitriangular left $H$-comodule algebra.
\end{proposition}

\begin{proof}
Since $K$ is invertible and $f$ is an algebra map, $K_C$ is invertible with inverse $(\id_H\otimes f)(K^{-1})$.
Applying $\id_H\otimes\id_H\otimes f$ to \eqref{eq:K-coproduct} and \eqref{eq:K-coaction} shows that $K_C$ satisfies the corresponding identities.
Now let $c\in C$.
We can choose $b\in B$ such that $f(b)=c$ by the surjectivity of $f$.
We have
\begin{align*}
K_C\delta_C(c)
&=(\id_H\otimes f)(K)
(\id_H\otimes f)(\delta_B(b))\\
&=(\id_H\otimes f)(K\delta_B(b))\\
&=(\id_H\otimes f)(\delta_B(b)K)
=\delta_C(c)K_C
\end{align*}
since $f$ is a left $H$-comodule algebra map.
Thus $K_C$ is a universal $K$-matrix on $C$.
\end{proof}

Quasitriangular $H$-comodule algebras form a category, defined as follows.

\begin{definition}[{\cite[Definition~6.14]{LWY25}}]
\label{def:QT-category}
The category ${}^H\mathsf{QT}$ of quasitriangular $H$-comodule algebras has quasitriangular $H$-comodule algebras as objects.
For two objects $(B,K),(B',K')\in{}^H\mathsf{QT}$, a morphism $\phi:(B,K)\to(B',K')$ is a left $H$-comodule algebra map $\phi:B\to B'$ satisfying
\begin{equation}
\label{eq:QT-map-K}
(\id_H\otimes\phi)(K)=K'.
\end{equation}
\end{definition}

\subsection{Cotensor products and Takeuchi's equivalence}
\label{subsec:quotients-cotensor}


Let $C$ be a coalgebra, $N$ a right $C$-comodule, and $M$ a left $C$-comodule.
The \emph{cotensor product} of $N$ and $M$ over $C$ is defined by
\begin{equation*}
N\square_C M
:=
\Ker(\rho_N\otimes\id_M-\id_N\otimes\delta_M)
\subseteq N\otimes M.
\end{equation*}
Let $H$ be a Hopf algebra and $L\subseteq H$ a \emph{left coideal subalgebra}, that is, a subalgebra of $H$ satisfying $\Delta(L)\subseteq H\otimes L$.
Set $L^+:=L\cap\Ker\varepsilon$.
Let $\pi:H\to H/HL^+$ be the canonical quotient coalgebra map onto the quotient left $H$-module coalgebra $H/HL^+$.
The map $(\id_H\otimes\pi)\circ\Delta:H\to H\otimes\pi(H)$ makes $H$ a right $\pi(H)$-comodule.
We define the space of right coinvariants of the right $\pi(H)$-comodule $H$ by
\begin{equation*}
	H^{\co\pi}
	:=
	\{h\in H\mid
	h_{(1)}\otimes\pi(h_{(2)})
	=
	h\otimes\pi(1_H)\}.
\end{equation*}

A \emph{relative $(H,L)$-Hopf module} is a left $H$-comodule $M$ equipped with a right $L$-module structure such that
\begin{equation*}
	\delta_M(ml)=m_{(-1)}l_{(1)}\otimes m_{(0)}l_{(2)}
	\quad(m\in M,\ l\in L).
\end{equation*}
For a relative $(H,L)$-Hopf module $M$, the quotient $M/ML^+$ is a left $\pi(H)$-comodule with coaction
\begin{equation*}
\overline{m}
\mapsto
\pi(m_{(-1)})\otimes\overline{m_{(0)}}.
\end{equation*}
Conversely, if $V$ is a left $\pi(H)$-comodule, then $H\square_{\pi(H)}V$ is a relative $(H,L)$-Hopf module with left $H$-coaction and right $L$-action induced by
\begin{equation*}
\delta_{H\square_{\pi(H)}V}(h\otimes v)
=
h_{(1)}\otimes(h_{(2)}\otimes v),\quad
(h\otimes v)\cdot\ell
=
h\ell\otimes v
\quad
(h\in H,\ \ell\in L,\ v\in V).
\end{equation*}

If $H$ is a finite-dimensional Hopf algebra, then its antipode is bijective.
Moreover, $H$ is free as both a left and a right module over any left coideal subalgebra $L\subseteq H$ by \cite[Theorem~6.1(ii)]{Skr07}.
Applying \cite[Theorem~1]{Tak79} to $H^{\cop}$, we obtain the following theorem.

\begin{theorem}[{\cite{Tak79,Skr07}}]
\label{thm:takeuchi-equivalence}
Let $H$ be a finite-dimensional Hopf algebra, $L\subseteq H$ a left coideal subalgebra, and $\pi:H\to\pi(H):=H/HL^+$ the canonical quotient coalgebra map.
Then the functors
\begin{equation*}
	M\mapsto M/ML^+,
	\quad
	V\mapsto H\square_{\pi(H)}V
\end{equation*}
give an adjoint equivalence between the category of relative $(H,L)$-Hopf modules and the category of left $\pi(H)$-comodules, with unit and counit given by the following natural isomorphisms, respectively:
\begin{gather*}
	M
	\to
	H\square_{\pi(H)}(M/ML^+),\quad
	m
	\mapsto
	m_{(-1)}\otimes\overline{m_{(0)}},
	\\
	\frac{H\square_{\pi(H)}V}{(H\square_{\pi(H)}V)L^+}
	\to
	V,\quad
	\overline{\sum_i h_i\otimes v_i}
	\mapsto
	\sum_i\varepsilon(h_i)v_i.
\end{gather*}
\end{theorem}


We obtain the following identity from \cref{thm:takeuchi-equivalence} by taking $M=L$ (see also \cite[Theorem~1.1(1)]{Tak94}):
\begin{equation}\label{eq:coideal-coinvariants}
	L=H^{\co\pi}.
\end{equation}

\subsection{$H$-simplicity, exactness, and indecomposability}
\label{subsec:comodule-module}

We use the standard terminology for finite tensor categories and their module categories; see \cite{EGNO15} for details.
A finite left module category $\mathcal M$ over a finite tensor category $\mathcal C$ is called \emph{exact} if $P\triangleright M$ is projective in $\mathcal M$ for every projective object $P\in\mathcal C$ and every object $M\in\mathcal M$.
It is called \emph{indecomposable} if it is nonzero and cannot be expressed as a direct sum of two nonzero $\mathcal C$-module subcategories.

Let $H$ be a finite-dimensional Hopf algebra and set $\mathcal C=H\text{-}\mathrm{mod}$.
An ideal $I\subseteq B$ of an $H$-comodule algebra $B$ is called \emph{$H$-costable} if $\delta_B(I)\subseteq H\otimes I$.
A nonzero left $H$-comodule algebra $B$ is called \emph{$H$-simple} if its only $H$-costable ideals are $0$ and $B$.
For a finite-dimensional $H$-comodule algebra $B$, the category $B\text{-}\mathrm{mod}$ is a module category over $\mathcal C$.
We call $B$ \emph{exact} if $B\text{-}\mathrm{mod}$ is exact as a $\mathcal C$-module category, and \emph{$H$-indecomposable} if this module category is indecomposable.
The latter condition is equivalent to the assertion that $B$ cannot be written as a direct sum of two nonzero $H$-costable ideals \cite[Lemma~3.9]{NSS25}.

We summarize the relation between $H$-simplicity, exactness, and indecomposability of $B$ in the following proposition.

\begin{proposition}
\label{prop:H-simple-exact-indecomposable}
Let $H$ be a finite-dimensional Hopf algebra and $B$ a finite-dimensional left $H$-comodule algebra.
Then the following conditions are equivalent:
\begin{enumerate}
\item[\rm{(i)}]
$B$ is $H$-simple.
\item[\rm{(ii)}]
$B$ is exact and $H$-indecomposable.
\end{enumerate}
\end{proposition}

\begin{proof}
Suppose that \textup{(i)} holds.
Then $B$ is a simple algebra object in ${}^{H}\mathrm{comod}$ and hence exact by the result recalled in \cite[Section~2.7]{NSS25}.
Thus $B$ is exact by \cite[Lemma~3.8]{NSS25} and $H$-indecomposable since it is $H$-simple.

Conversely, suppose that \textup{(ii)} holds.
The algebra $B$ is an exact algebra object in ${}^{H}\mathrm{comod}$ by the same lemma.
It follows from \cite[Theorem~7.1]{CSZ25} that there exist simple algebra objects $B_1,\ldots,B_n$ in ${}^{H}\mathrm{comod}$ such that
$B\cong B_1\times\cdots\times B_n$ as algebra objects.
Since $B$ is $H$-indecomposable, we have $n=1$.
Thus $B$ is a simple algebra object in ${}^{H}\mathrm{comod}$.
The ideals of $B$ as an algebra object in ${}^{H}\mathrm{comod}$ are precisely its $H$-costable ideals, so $B$ is $H$-simple.
\end{proof}

\section{Reflective algebras and factorizability}
\label{sec:reflective-centers-algebras}

Laugwitz, Walton and Yakimov \cite{LWY25} introduced the reflective center $\mathcal E_{\mathcal C}(\mathcal M)$ for a braided monoidal category $\mathcal C$ and a left $\mathcal C$-module category $\mathcal M$.
It is an analogue of the Drinfeld center and carries a canonical braided $\mathcal C$-module category structure.
Since we do not use the categorical definition in this paper, we restrict our attention to results concerning comodule algebras in this section.
We refer to \cite[Section~4]{LWY25} for the detailed definition of the reflective center and its categorical structure.

In \cite{LWY25}, the reflective algebra $R_H(A)$ is constructed for a quasitriangular Hopf algebra $(H,R)$ and a left $H$-comodule algebra $A$ as the Hopf algebra counterpart of the reflective center.
In this section, we recall the basic properties of reflective algebras and their factorizability.
We first describe the construction of reflective algebras and their canonical quasitriangular $H$-comodule algebra structures.
We then present the universal property of $R_H(\bfk)$ and its relation to the Drinfeld map.
Finally, we recall the definition of factorizable quasitriangular comodule algebras and the relevant results of Walton and Yadav \cite{WY25}.

\subsection{The definition of reflective algebras}
\label{subsec:reflective-algebras}

Throughout this section, let $(H,R)$ be a quasitriangular Hopf algebra and $A$ a left $H$-comodule algebra.
Write $R=\sum_i s_i\otimes t_i\in H\otimes H$.

Let $\widehat H$ be the vector space $H$ equipped with the following comultiplication, counit, and left $H$-action formulas, for $h\in\widehat H$ and $\ell\in H$:
\begin{align*}
\widehat\Delta(h) = \sum_{i,j}t_jh_{(1)}t_i\otimes h_{(2)}s_iS^{-1}(s_j),\quad
\widehat\varepsilon(h) = \varepsilon(h),\quad
\ell\rightharpoonup h = \ell_{(2)}hS^{-1}(\ell_{(1)}).
\end{align*}
These maps make $\widehat H$ a left $H$-module coalgebra.
Since $H$ is finite-dimensional, the transpose of $\widehat\Delta$ defines an algebra structure on $\widehat H^*$.
The action given by
\begin{equation}\label{eq:hatHstar-right-H-action}
\langle\xi\leftharpoonup\ell,h\rangle
=
\langle\xi,\ell\rightharpoonup h\rangle
\quad
(\xi\in\widehat H^*,\ h\in\widehat H,\ \ell\in H)
\end{equation}
makes $\widehat H^*$ a right $H$-module algebra.
Consequently, $(\widehat H^*)^{\op}$ is a right $H^{\cop}$-module algebra with the same action.

\begin{definition}[{\cite[Definition~5.27]{LWY25}}]
\label{def:reflective-algebra}
The \emph{reflective algebra} of $A$ with respect to $H$ is the crossed product algebra
\begin{equation*}
R_H(A):=
A\rtimes_H(\widehat H^*)^{\op}.
\end{equation*}
As vector spaces, $R_H(A)\cong A\otimes(\widehat H^*)^{\op}$.
Under this identification, the multiplication is given by
\begin{equation}\label{eq:reflective-algebra-product}
(a\xi)(b\zeta)
=
a b_{(0)}\,(\xi\leftharpoonup b_{(-1)})\zeta
\quad
(a,b\in A,\ \xi,\zeta\in(\widehat H^*)^{\op}).
\end{equation}
Here $a\,\xi$ denotes $a\otimes\xi\in A\otimes(\widehat H^*)^{\op}$.
\end{definition}

By the canonical embeddings into the crossed product, we may regard $A$ and $(\widehat H^*)^{\op}$ as subalgebras of $R_H(A)$.
Since $\widehat H^*$ and $H^*$ are the same vector space, we also have an isomorphism $R_H(A)\cong A\otimes H^*$ of vector spaces.
In particular, if $A=\bfk$ is equipped with the trivial left $H$-coaction, then $R_H(\bfk)\cong(\widehat H^*)^{\op}$ as $\bfk$-algebras.

We next recall the canonical quasitriangular structure on a reflective algebra and its relation to the reflective center.
Let $\{h_d,\xi_d\}_d$ be a pair of dual bases of $H$ and $H^*$, respectively, so that $h_d\in H$, $\xi_d\in H^*$, and
\begin{equation}\label{eq:dual-basis}
h=\sum_d\langle\xi_d,h\rangle h_d,
\quad
\xi=\sum_d\langle\xi,h_d\rangle\xi_d
\end{equation}
for all $h\in H$ and $\xi\in H^*$.

\begin{theorem}[{\cite[Corollary~6.9]{LWY25}}]
\label{thm:reflective-algebra-structure}
The algebra $R_H(A)$ has a left $H$-comodule algebra structure
\begin{equation*}
\delta_{\mathrm{ref}}:
R_H(A)\to H\otimes R_H(A)
\end{equation*}
given, for $a\in A$ and $\xi\in(\widehat H^*)^{\op}$, by
\begin{align*}
\delta_{\mathrm{ref}}(a) = a_{(-1)}\otimes a_{(0)},\hspace{0.8em}
\delta_{\mathrm{ref}}(\xi) = \sum_{i,j,d}\langle\xi,t_jh_ds_i\rangle\,s_jt_i\otimes\xi_d,\hspace{0.8em}
\delta_{\mathrm{ref}}(a\xi) = \delta_{\mathrm{ref}}(a)\delta_{\mathrm{ref}}(\xi).
\end{align*}
Moreover,
\begin{equation}\label{eq:K-ref}
K_{\mathrm{ref}}(A)
:=
\sum_d h_d\otimes\xi_d
\in
H\otimes(\widehat H^*)^{\op}
\subseteq
H\otimes R_H(A)
\end{equation}
is a universal $K$-matrix for $R_H(A)$.
Thus $(R_H(A),\delta_{\mathrm{ref}},K_{\mathrm{ref}}(A))$ is a quasitriangular left $H$-comodule algebra.
\end{theorem}

The coaction $\delta_{\mathrm{ref}}$ makes $R_H(A)\text{-}\mathrm{mod}$ a module category over $H\text{-}\mathrm{mod}$.
The universal $K$-matrix $K_{\mathrm{ref}}(A)$ equips it with a braided module category structure.

\subsection{The universal property of \texorpdfstring{$R_H(\bfk)$}{RH(k)} and the Drinfeld map}
\label{subsec:universal-property-drinfeld-map}

The reflective algebra $(R_H(\bfk),K_{\mathrm{ref}}(\bfk))$ is an initial object in the category ${}^H\mathsf{QT}$ of quasitriangular left $H$-comodule algebras \cite[Theorem~6.15]{LWY25}.
For a quasitriangular left $H$-comodule algebra $(B,K)$, we denote by
\begin{equation}\label{eq:QT-initial-map}
\kappa_{(B,K)}:
(R_H(\bfk),K_{\mathrm{ref}}(\bfk))
\to
(B,K)
\end{equation}
the unique morphism in ${}^H\mathsf{QT}$ from $(R_H(\bfk),K_{\mathrm{ref}}(\bfk))$ to $(B,K)$.

Regarding $H$ as a left $H$-comodule algebra via its coproduct $\Delta:H\to H\otimes H$, the element $R_{21}R\in H\otimes H$ is a universal $K$-matrix for $H$.
Thus $(H,R_{21}R)$ is an object of ${}^H\mathsf{QT}$ \cite[Example~6.12(b)]{WY25}.
Applying the universal property of $R_H(\bfk)$ to this object gives the following proposition.

\begin{proposition}
\label{prop:drinfeld-map-QT-map}
For the object $(H,R_{21}R)$ of ${}^H\mathsf{QT}$, the unique quasitriangular left $H$-comodule algebra map
$\kappa_{(H,R_{21}R)}:(R_H(\bfk),K_{\mathrm{ref}}(\bfk))\to(H,R_{21}R)$
induced by the universal property of $R_H(\bfk)$ coincides with the Drinfeld map $Q_H$ under the identification $R_H(\bfk)\cong H^*$ of vector spaces.
In particular, $Q_H$ is a quasitriangular left $H$-comodule algebra map.
\end{proposition}

\begin{proof}
Let $\{h_d,\xi_d\}_d$ be a pair of dual bases of $H$ and $H^*$.
By \eqref{eq:K-ref} and \eqref{eq:QT-map-K}, we have
$\sum_d h_d\otimes\kappa_{(H,R_{21}R)}(\xi_d)=R_{21}R$.
Hence, by \eqref{eq:dual-basis}, for every $\xi\in H^*$ we obtain
\begin{equation*}
\kappa_{(H,R_{21}R)}(\xi)
=
(\xi\otimes\id_H)(R_{21}R)
=
Q_H(\xi).
\end{equation*}
Therefore $\kappa_{(H,R_{21}R)}=Q_H$, as required.
\end{proof}

\begin{corollary}
\label{cor:factorizable-H-drinfeld-map}
The following conditions are equivalent.
\begin{enumerate}
\item[\rm{(i)}] The quasitriangular Hopf algebra $(H,R)$ is factorizable.
\item[\rm{(ii)}] The map $Q_H$ is an isomorphism from $(R_H(\bfk),K_{\mathrm{ref}}(\bfk))$ to $(H,R_{21}R)$ in ${}^H\mathsf{QT}$.
\item[\rm{(iii)}] The object $(H,R_{21}R)$ is also an initial object of ${}^H\mathsf{QT}$.
\end{enumerate}
\end{corollary}

\begin{proof}
This follows from the definition of factorizability, \cref{prop:drinfeld-map-QT-map}, and the uniqueness of initial objects up to isomorphism.
\end{proof}

\subsection{Factorizability of quasitriangular comodule algebras}
\label{subsec:factorizability-quasitriangular-comodule-algebras}

Walton and Yadav introduced factorizability for quasitriangular left $H$-comodule algebras and studied its properties \cite{WY25}.
In this subsection, we recall the definitions and results needed later.
For the categorical definitions concerning braided module categories and their relation to this notion, see \cite[Section~5 and Theorem~6.19]{WY25}.

Let $(B,K)$ be a quasitriangular left $H$-comodule algebra with $K=\sum_i K_i\otimes K^i$.
Define the subspace $E(H,B)$ of $\Hom_{\bfk}(H,B)$ by
\begin{equation} \label{eq:E-H-B}
E(H,B)
:=
\{\xi\in\Hom_{\bfk}(H,B)\mid
\xi(b_{(-1)}h)b_{(0)}=b\xi(h)
\text{ for all }b\in B,\ h\in H\}.
\end{equation}
Furthermore, define the linear map $\theta_B:H^*\to E(H,B)$ by
\begin{equation} \label{eq:theta-B}
\theta_B(f)(h)
=
\sum_i\langle f,S(h_{(1)})K_i h_{(2)}\rangle K^i
\quad(f\in H^*,\ h\in H).
\end{equation}
The map $\theta_B$ is well-defined by \cite[Lemma~6.16]{WY25}.
Following Walton and Yadav, we define a factorizable comodule algebra as follows.

\begin{definition}[{\cite[Hypothesis~6.14 and Definition~6.15]{WY25}}]
\label{def:factorizable-comodule-algebra}
Let $(B,K)$ be an exact and $H$-indecomposable quasitriangular left $H$-comodule algebra.
The pair $(B,K)$ is called a \emph{factorizable left $H$-comodule algebra} if the linear map $\theta_B:H^*\to E(H,B)$ is an isomorphism.
\end{definition}

The subspace $E(H,B)$ and the linear map $\theta_B$ can be defined without assuming that $B$ is exact or $H$-indecomposable.
Under these assumptions, however, the following holds.

\begin{lemma}[{\cite[Lemma~6.17]{WY25}}]
\label{lem:factorizable-comodule-injective}
If $B$ is exact and $H$-indecomposable, then $\dim E(H,B)=\dim H$.
Consequently, $(B,K)$ is factorizable if and only if $\theta_B:H^*\to E(H,B)$ is injective.
\end{lemma}

For reflective algebras, the following additional results are known.
Let $(H,R)$ be a quasitriangular Hopf algebra and $A$ a left $H$-comodule algebra.
\begin{itemize}
\item The map $\theta_{R_H(A)}:H^*\to E(H,R_H(A))$ is injective. This does not require $R_H(A)$ to be $H$-simple \cite[Proof of Proposition~6.24]{WY25}.
\item If $R_H(A)$ is $H$-simple, then $(R_H(A),K_{\mathrm{ref}}(A))$ is a factorizable left $H$-comodule algebra \cite[Proposition~6.24]{WY25}.
\end{itemize}

By \cref{prop:H-simple-exact-indecomposable}, the exactness and $H$-indecomposability of $R_H(A)$ are equivalent to its $H$-simplicity.
Therefore the above results yield the following corollary.

\begin{corollary}
\label{cor:reflective-factorizable-H-simple}
The following conditions are equivalent.
\begin{enumerate}
\item[\rm{(i)}] $R_H(A)$ is $H$-simple.
\item[\rm{(ii)}] $(R_H(A),K_{\mathrm{ref}}(A))$ is a factorizable left $H$-comodule algebra.
\end{enumerate}
\end{corollary}

\begin{proof}
The implication \textup{(i)} $\Rightarrow$ \textup{(ii)} follows from \cite[Proposition~6.24]{WY25}.
Conversely, suppose that \textup{(ii)} holds.
Then $R_H(A)$ is exact and $H$-indecomposable by the definition of a factorizable comodule algebra.
Thus \cref{prop:H-simple-exact-indecomposable} implies \textup{(i)}.
\end{proof}

\section{Quotients of reflective algebras associated with the Drinfeld map}
\label{sec:drinfeld-quotient}

In this section, let $(H,R)$ be a quasitriangular Hopf algebra and $A$ a left $H$-comodule algebra.
We use the Drinfeld map $Q_H$ defined in \cref{def:drinfeld-map-factorizable}.
By \cref{prop:drinfeld-map-QT-map}, the map $Q_H:R_H(\bfk)\to H$ is a left $H$-comodule algebra map.
We realize the quotient $R_H(A)/A\,\Ker Q_H$ as a cotensor product.
We also determine when this quotient is $H$-simple.

\subsection{The quotient associated with the Drinfeld map}
\label{subsec:drinfeld-quotient}

As in \cref{sec:reflective-centers-algebras}, we identify $R_H(\bfk)$ with $H^*$ as vector spaces.
We also regard $A$ and $R_H(\bfk)$ as subalgebras of $R_H(A)$.
Since $Q_H$ is a left $H$-comodule algebra map, its kernel $\Ker Q_H$ is an ideal and a left $H$-subcomodule of $R_H(\bfk)$.
The image $\Img Q_H$ is a left coideal subalgebra of $H$.
A left coideal subalgebra $L$ of $H$ is called normal if $h_{(1)}xS(h_{(2)})\in L$ for all $h\in H$ and $x\in L$.

\begin{lemma}
\label{lem:kernel-image-QH}
The following statements hold.
\begin{enumerate}
\item[\rm{(i)}] $(\Ker Q_H)\leftharpoonup H\subseteq\Ker Q_H$.
\item[\rm{(ii)}] The subalgebra $\Img Q_H$ is a normal left coideal subalgebra of $H$.
\end{enumerate}
\end{lemma}

\begin{proof}
By \cite[Proposition~2.1.14]{Maj95} and \eqref{eq:hatHstar-right-H-action}, for all $\xi\in R_H(\bfk)$ and $h\in H$, we have
\begin{equation}\label{eq:QH-right-H-action}
Q_H(\xi\leftharpoonup h)
=
S^{-1}(h_{(2)})Q_H(\xi)h_{(1)}.
\end{equation}
Consequently, if $\xi\in\Ker Q_H$, then $Q_H(\xi\leftharpoonup h)=0$, which proves \textup{(i)}.
Next, let $x=Q_H(\xi)\in\Img Q_H$ and $h\in H$. Replacing $h$ with $S(h)$ in \eqref{eq:QH-right-H-action}, we obtain
\begin{equation*}
h_{(1)}xS(h_{(2)})
=
Q_H(\xi\leftharpoonup S(h))
\in
\Img Q_H.
\end{equation*}
Thus $\Img Q_H$ is stable under the left adjoint action of $H$, and hence is a normal left coideal subalgebra of $H$.
\end{proof}

\begin{proposition}
\label{prop:canonical-ideal}
The subspace $A\,\Ker Q_H\subseteq R_H(A)$ is an $H$-costable ideal.
\end{proposition}

\begin{proof}
Since $\Ker Q_H$ is an ideal of $R_H(\bfk)$, the subspace $A\,\Ker Q_H$ is closed under left and right multiplication by elements of $R_H(\bfk)$ and under left multiplication by elements of $A$.
Moreover, for all $a,b\in A$ and $\xi\in\Ker Q_H$, \eqref{eq:reflective-algebra-product} and \cref{lem:kernel-image-QH}(i) give
\begin{equation*}
(a\xi)b
=
a b_{(0)}(\xi\leftharpoonup b_{(-1)})
\in
A\,\Ker Q_H.
\end{equation*}
Since $A$ and $R_H(\bfk)$ generate $R_H(A)$ as an algebra, it follows that $A\,\Ker Q_H$ is an ideal of $R_H(A)$.
Furthermore, $A$ and $\Ker Q_H$ are left $H$-subcomodules of $R_H(A)$, and the coaction on $R_H(A)$ is an algebra map. Therefore,
$ \delta_{\mathrm{ref}}(A\,\Ker Q_H)
\subseteq
H\otimes A\,\Ker Q_H $
and hence $A\,\Ker Q_H$ is $H$-costable.
\end{proof}

By \cref{prop:canonical-ideal}, the quotient $R_H(A)/A\,\Ker Q_H$ is a left $H$-comodule algebra.
The following proposition shows that this quotient has a quasitriangular structure.
\begin{proposition}
	\label{prop:quotient-cotensor-K-matrices}
	Let $\{h_d,\xi_d\}_d$ be a pair of dual bases of $H$ and $H^*$, and let $\overline{\xi_d}$ denote the image of $\xi_d\in R_H(\bfk)\subset R_H(A)$ under the canonical quotient map $R_H(A)\to R_H(A)/A\,\Ker Q_H$.
	Set
	\begin{equation*}
		\overline{K}_{\mathrm{ref}}(A) :=
		\sum_d h_d\otimes\overline{\xi_d}
		\in
		H\otimes(R_H(A)/A\,\Ker Q_H).
	\end{equation*}
	Then $(R_H(A)/A\,\Ker Q_H,\overline{K}_{\mathrm{ref}}(A))$ is a quasitriangular left $H$-comodule algebra.
\end{proposition}

\begin{proof}
	The ideal $A\,\Ker Q_H$ is $H$-costable.
	Thus the quotient map is a surjective left $H$-comodule algebra map.
	Since $K_{\mathrm{ref}}(A)=\sum_d h_d\otimes\xi_d$ is a universal $K$-matrix for $R_H(A)$, \cref{prop:K-matrix-quotient} implies that $\overline{K}_{\mathrm{ref}}(A)$ is a universal $K$-matrix for $R_H(A)/A\,\Ker Q_H$.
\end{proof}

The subspace $A\,\Ker Q_H\subset R_H(A)$ is identified with $A\otimes\Ker Q_H$ as a vector space. Thus, we obtain an isomorphism of vector spaces
\begin{equation*}
	R_H(A)/A\,\Ker Q_H
	\cong
	A\otimes\Img Q_H.
\end{equation*}

\subsection{The quotient of $R_H(A)$ as a cotensor product}
\label{subsec:cotensor-description}

Set $(\Img Q_H)^+:=\Img Q_H\cap\Ker\varepsilon$.
By \cref{lem:kernel-image-QH}(ii), the subalgebra $\Img Q_H$ is a normal left coideal subalgebra of $H$.
By \cite[Lemma~4.2]{Tak94}, we have
$ H(\Img Q_H)^+
=
(\Img Q_H)^+H $
and this is a Hopf ideal of $H$.
Let the quotient Hopf algebra $T$ be defined by
\begin{equation}\label{eq:T-def}
T:=H/H(\Img Q_H)^+
\end{equation}
and denote the quotient Hopf algebra map by $\pi:H\to T$.
Define a right $T$-coaction on $H$ and a left $T$-coaction on $A$, respectively, by
\begin{align*}
\delta_H^T:H\to H\otimes T,\quad
&\delta_H^T(h)=(\id_H\otimes\pi)\Delta(h),\\
\delta_A^T:A\to T\otimes A,\quad
&\delta_A^T(a)=(\pi\otimes\id_A)\delta_A(a).
\end{align*}
Then the cotensor product is defined by
\begin{equation*}
H\square_T A
:=
\{z\in H\otimes A\mid
(\delta_H^T\otimes\id_A)(z)
=
(\id_H\otimes\delta_A^T)(z)\}.
\end{equation*}
Since $\delta_H^T$ and $\delta_A^T$ are algebra maps, $H\square_T A$ is a subalgebra of $H\otimes A$.
Moreover, the restriction of $\Delta\otimes\id_A$ to $H\square_T A$ gives a natural left $H$-coaction, making $H\square_T A$ a left $H$-comodule algebra, as follows:
\begin{equation*}
\delta_{H\square_T A}(h\otimes a)
=
h_{(1)}\otimes(h_{(2)}\otimes a).
\end{equation*}

The main result of this section is the following.
\begin{theorem}
	\label{thm:cotensor-quasitriangular}
	The map
	\begin{equation*}
		\phi:
		R_H(A)/A\,\Ker Q_H
		\to
		H\square_T A,
		\quad
		\phi(\overline{a\xi})
		=
		a_{(-1)}Q_H(\xi)\otimes a_{(0)}
	\end{equation*}
	induces the following isomorphism of quasitriangular left $H$-comodule algebras:
	\begin{equation*}
		(R_H(A)/A\,\Ker Q_H, \overline{K}_{\mathrm{ref}}(A))
		\cong (H\square_T A,(R_{21}R)_{12}).
	\end{equation*}
\end{theorem}
\begin{proof}
	We first show that $\phi$ is an isomorphism of left $H$-comodules using \cref{thm:takeuchi-equivalence}.
	Set $L=\Img Q_H$.
	By \cref{lem:kernel-image-QH}(ii), $L$ is a left coideal subalgebra of $H$, and we have $T=H/HL^+$.
	Equip $A\otimes L$ with the left $H$-coaction and the right $L$-action given by
	\begin{equation*}
		\delta(a\otimes x)
		=a_{(-1)}x_{(1)}\otimes(a_{(0)}\otimes x_{(2)}),\quad
		(a\otimes x)\cdot y
		=a\otimes xy
	\end{equation*}
	for $a\in A$ and $x,y\in L$.
	Then $A\otimes L$ is a relative $(H,L)$-Hopf module.
	Since $\pi(hx)=\varepsilon(x)\pi(h)$ for all $h\in H$ and $x\in L$, the map $\overline{a\otimes x}\mapsto\varepsilon(x)a$ identifies $(A\otimes L)/(A\otimes L)L^+$ with $A$ as left $T$-comodules, where $A$ carries the coaction $\delta_A^T$.
	Under this identification, the unit of \cref{thm:takeuchi-equivalence} becomes
	\begin{equation*}
		A\otimes L\to H\square_T A,\quad
		a\otimes x\mapsto a_{(-1)}x\otimes a_{(0)}.
	\end{equation*}
	Thus this map is an isomorphism of left $H$-comodules.
	Moreover, the map $\overline{a\xi}\mapsto a\otimes Q_H(\xi)$ identifies $R_H(A)/A\,\Ker Q_H$ with $A\otimes L$ as left $H$-comodules, since $Q_H$ is a left $H$-comodule map and $\delta_{\mathrm{ref}}$ is an algebra map.
	Consequently, $\phi$ is an isomorphism of left $H$-comodules.

	We next show that $\phi$ is an algebra map.
	Let $a,b\in A$ and $\xi,\zeta\in R_H(\bfk)$.
	By \eqref{eq:reflective-algebra-product} and \eqref{eq:QH-right-H-action}, we have
	\begin{align*}
		\phi(\overline{a\xi}\,\overline{b\zeta})
		&=
		\phi(\overline{ab_{(0)}(\xi\leftharpoonup b_{(-1)})\zeta})\\
		&=
		a_{(-1)}b_{(-1)}
		Q_H(\xi\leftharpoonup b_{(-2)})Q_H(\zeta)
		\otimes a_{(0)}b_{(0)}\\
		&=
		a_{(-1)}b_{(-1)}S^{-1}(b_{(-2)})
		Q_H(\xi)b_{(-3)}Q_H(\zeta)
		\otimes a_{(0)}b_{(0)}\\
		&=
		a_{(-1)}Q_H(\xi)b_{(-1)}Q_H(\zeta)
		\otimes a_{(0)}b_{(0)}\\
		&=
		\phi(\overline{a\xi})\phi(\overline{b\zeta}).
	\end{align*}
	Hence $\phi$ is an isomorphism of left $H$-comodule algebras.

	Finally, we show that $\phi$ preserves the universal $K$-matrices.
	Let $\{h_d,\xi_d\}_d$ be a pair of dual bases of $H$ and $H^*$.
	By the definition of the Drinfeld map, we have
	\begin{equation*}
		(\id_H\otimes\phi)(\overline{K}_{\mathrm{ref}}(A))
		=
		\sum_d h_d\otimes\phi(\overline{\xi_d})
		=
		\sum_d h_d\otimes Q_H(\xi_d)\otimes1_A
		=
		(R_{21}R)_{12}.
	\end{equation*}
	Since $\overline{K}_{\mathrm{ref}}(A)$ is a universal $K$-matrix by \cref{prop:quotient-cotensor-K-matrices} and $\phi$ is an isomorphism of left $H$-comodule algebras, the element $(R_{21}R)_{12}$ is a universal $K$-matrix for $H\square_T A$.
	Therefore, $\phi$ is an isomorphism of quasitriangular left $H$-comodule algebras.
\end{proof}

\begin{remark}
	\label{rem:muger-center-quotient}
	The quotient Hopf algebra map $\pi:H\to T$ induces a fully faithful functor from $T\text{-}\mathrm{mod}$ to $H\text{-}\mathrm{mod}$.
	Its essential image is the M\"{u}ger center of $H\text{-}\mathrm{mod}$ \cite[Lemma~3.3]{WZ25}.
	For the semisimple case over an algebraically closed field of characteristic zero, see also \cite[Proposition~4.5(1)]{Bur20}.
\end{remark}

\subsection{Costable ideals and $H$-simplicity of the quotient}
\label{subsec:costable-ideals-simplicity}
In this subsection, we use $H\square_T A$ to characterize when the quotient $R_H(A)/A\,\Ker Q_H$ studied in \cref{subsec:cotensor-description} is $H$-simple.

\begin{theorem}
\label{thm:costable-ideal-correspondence}
Let $T$ be defined by \eqref{eq:T-def}.
Then there is a bijective correspondence between the $T$-costable ideals of $A$ and the $H$-costable ideals of $H\square_T A$.
\end{theorem}

\begin{proof}
	By \cref{lem:kernel-image-QH}(ii), $\Img Q_H$ is a left coideal subalgebra of $H$.
	By \eqref{eq:T-def}, we also have $T=H/H(\Img Q_H)^+$.
	Thus we may apply \cref{thm:takeuchi-equivalence} with $L=\Img Q_H$ and $\pi(H)=T$.

	Let $I\subset A$ be a $T$-costable ideal.
	Since $I$ is a left $T$-comodule, \cref{thm:takeuchi-equivalence} associates it with
	\begin{equation*}
		H\square_T I\subset H\square_TA.
	\end{equation*}
	This is a relative $(H,\Img Q_H)$-Hopf module and, in particular, a left $H$-subcomodule of $H\square_T A$.
	Since $H\square_T A$ is a subalgebra of $H\otimes A$ and $I$ is an ideal of $A$, the subspace $H\square_T I$ is an ideal of $H\square_T A$.
	Hence $H\square_T I$ is an $H$-costable ideal of $H\square_T A$.

	Conversely, let $J\subset H\square_T A$ be an $H$-costable ideal.
	Applying \eqref{eq:coideal-coinvariants} with $L=\Img Q_H$ and $\pi(H)=T$, we have $\Img Q_H=H^{\co\pi}$.
	Thus, for any $x\in\Img Q_H$, we have $x\otimes1_A\in H\square_T A$, and since $J$ is an ideal, $J(x\otimes1_A)\subset J$.
	Consequently, restricting the relative $(H,\Img Q_H)$-Hopf module structure of $H\square_T A$ to $J$ makes $J$ a relative $(H,\Img Q_H)$-Hopf module.
	By the natural isomorphism in \cref{thm:takeuchi-equivalence}, the left $T$-comodule corresponding to $J$ is identified with
	\begin{equation*}
		(\varepsilon\otimes\id_A)(J)\subset A.
	\end{equation*}
	Furthermore, the algebra map $\varepsilon\otimes\id_A:H\square_T A\to A$ is surjective.
	Indeed, for every $a\in A$, we have $\delta_A(a)\in H\square_T A$ and $(\varepsilon\otimes\id_A)(\delta_A(a))=a$.
	Therefore, $(\varepsilon\otimes\id_A)(J)$ is an ideal of $A$.
	As noted above, it is also a left $T$-comodule and hence a $T$-costable ideal of $A$.
	Thus, the correspondence of \cref{thm:takeuchi-equivalence} restricts to a bijective correspondence between the $T$-costable ideals of $A$ and the $H$-costable ideals of $H\square_T A$.

\end{proof}

\begin{corollary}
\label{cor:quotient-H-simple}
The following conditions are equivalent.
\begin{enumerate}
\item[\rm{(i)}] The quotient $R_H(A)/A\,\Ker Q_H$ is $H$-simple.
\item[\rm{(ii)}] The algebra $A$ is $T$-simple.
\end{enumerate}
\end{corollary}

\begin{proof}
This follows from \cref{thm:cotensor-quasitriangular,thm:costable-ideal-correspondence}.
\end{proof}

\section{Quasitriangular comodule algebras over factorizable Hopf algebras}
\label{sec:factorizable-case}

Throughout this section, let $(H,R)$ be a factorizable Hopf algebra and $(B,K)$ a nonzero quasitriangular left $H$-comodule algebra.
For algebras $A,A'$ and left $H$-comodule algebras $C,C'$, let $\mathrm{Alg}(A,A')$ and $\HAlg(C,C')$ denote the sets of algebra maps from $A$ to $A'$ and left $H$-comodule algebra maps from $C$ to $C'$, respectively.
The center of an algebra $A$ is denoted by $Z(A)$.

\subsection{Structure of quasitriangular comodule algebras}
\label{subsec:comodalg-over-factorizable-H}

By the universal property of $R_H(\bfk)$ recalled in \cref{subsec:universal-property-drinfeld-map}, there exists a unique quasitriangular left $H$-comodule algebra map
$ \kappa_{(B,K)}$
as in \eqref{eq:QT-initial-map}.
Moreover, by \cref{prop:drinfeld-map-QT-map}, the Drinfeld map
$Q_H
$
is an isomorphism of quasitriangular left $H$-comodule algebras.
Set
\begin{equation}\label{eq:jB}
j_B
:=
\kappa_{(B,K)}\circ Q_H^{-1}
:
H\to B.
\end{equation}
Thus $j_B:H\to B$ is a left $H$-comodule algebra map, and the following identities hold:
\begin{align}\label{eq:jB-colinear}
\delta_B(j_B(h))
&=
h_{(1)}\otimes j_B(h_{(2)}),
\\
K
&=
(\id_H\otimes j_B)(R_{21}R)
\label{eq:KB-jB}
\end{align}
for every $h\in H$.

\begin{lemma}
\label{lem:jB-commutes-coinvariants}
With the notation above, for all $c\in B^{\co H}$ and $h\in H$, we have $ j_B(h)c=cj_B(h) $.
\end{lemma}

\begin{proof}
Let $c$ be an arbitrary element of $B^{\co H}$.
By the defining relation \eqref{eq:K-commutative} for $K$, we have $ K(1_H\otimes c)=(1_H\otimes c)K $.
Applying $f\otimes\id_B$ to this equality for any $f\in H^*$, we obtain
\begin{equation*}
	(f\otimes\id_B)(K)c
	=
	c(f\otimes\id_B)(K).
\end{equation*}
Equation \eqref{eq:KB-jB} and the definition of the Drinfeld map also give
\begin{equation}\label{eq:KB-slice-QH}
	(f\otimes\id_B)(K)
	=
	j_B((f\otimes\id_H)(R_{21}R))
	=
	j_B(Q_H(f)).
\end{equation}
Hence $ j_B(Q_H(f))c=cj_B(Q_H(f)) $ for all $f\in H^*$. Since $(H,R)$ is factorizable, $Q_H:H^*\to H$ is surjective.
This completes the proof.
\end{proof}

We define a left $H$-coaction on $B^{\co H}\otimes H$ as follows:
\begin{equation}\label{eq:coinvariants-tensor-H-coaction}
\delta_{B^{\co H}\otimes H} (c\otimes h)
=
h_{(1)}\otimes(c\otimes h_{(2)})
\quad
(c\in B^{\co H},\ h\in H).
\end{equation}
It is straightforward to verify that
$ (R_{21}R)_{13}
\in
H\otimes B^{\co H}\otimes H $
is a universal $K$-matrix for $B^{\co H}\otimes H$.

\begin{theorem}\label{thm:factorizable-B-structure}
	Let $(H,R)$ be a factorizable Hopf algebra and $(B,K)$ a quasitriangular left $H$-comodule algebra.
	The map
	\begin{equation*}
		F_B:
		B^{\co H}\otimes H
		\to
		B,
		\quad
		c\otimes h \mapsto c\,j_B(h)
	\end{equation*}
	induces the following isomorphism of quasitriangular left $H$-comodule algebras:
	\begin{equation*}
		(B^{\co H}\otimes H,(R_{21}R)_{13})
		\cong (B,K).
	\end{equation*}

\end{theorem}

\begin{proof}
First, we have
\begin{equation*}
	F_B((c\otimes h)(d\otimes k))
	=
	cd\,j_B(hk)
	=
	c\,j_B(h)d\,j_B(k)
	=
	F_B(c\otimes h)F_B(d\otimes k)
\end{equation*}
for all $c,d\in B^{\co H}$ and $h,k\in H$ by \cref{lem:jB-commutes-coinvariants}.
Hence $F_B$ is an algebra map.
Next, define a right $H$-action on $B$ by
\begin{equation}\label{eq:right-action-on-B}
b\triangleleft h
:=
b\,j_B(h)
\quad
(b\in B,\ h\in H).
\end{equation}
By \eqref{eq:jB-colinear}, we have
\begin{equation*}
\delta_B(b\triangleleft h)
=
b_{(-1)}h_{(1)}
\otimes
(b_{(0)}\triangleleft h_{(2)}).
\end{equation*}
Thus the right $H$-action and the left $H$-coaction make $B$ a Hopf module. By the Fundamental Theorem of Hopf Modules (see, for example, \cite[Theorem~1.9.4]{Mon93}), $F_B$ is an isomorphism of left $H$-comodules. Since $F_B$ is also an algebra map, it is an isomorphism of left $H$-comodule algebras. Finally, \eqref{eq:KB-jB} gives
\begin{equation*}
(\id_H\otimes F_B)((R_{21}R)_{13})
=
(\id_H\otimes j_B)(R_{21}R)
=
K.
\end{equation*}
Therefore, $F_B$ is an isomorphism of quasitriangular left $H$-comodule algebras.
\end{proof}

\cref{thm:factorizable-B-structure} also characterizes when a universal $K$-matrix exists on a left $H$-comodule algebra.

\begin{corollary}\label{cor:K-matrix-existence}
	Let $(H,R)$ be a finite-dimensional factorizable Hopf algebra and $B$ a finite-dimensional left $H$-comodule algebra.
	Then there exists a universal $K$-matrix on $B$ if and only if $B$ is isomorphic to $B^{\co H}\otimes H$ as a left $H$-comodule algebra.
	In this case, $\dim B=\dim B^{\co H}\cdot\dim H$, and hence $\dim H$ divides $\dim B$.
\end{corollary}

\begin{proof}
	The sufficiency follows from \cref{thm:factorizable-B-structure}.
	Conversely, $(R_{21}R)_{13}$ is a universal $K$-matrix on $B^{\co H}\otimes H$.
	An isomorphism of left $H$-comodule algebras from $B^{\co H}\otimes H$ to $B$ carries this universal $K$-matrix to a universal $K$-matrix on $B$, proving the necessity.
\end{proof}

\subsection{Classification of universal $K$-matrices}
\label{subsec:comodalg-over-factorizable-H-Kmatrix}

In this subsection, we prove the following classification of universal $K$-matrices on the left $H$-comodule algebra $B$.

\begin{theorem}\label{thm:K-matrices-general-B}
Let $(B,K)$ be a quasitriangular left $H$-comodule algebra.
Then the universal $K$-matrices on $B$ are in bijection with
$ \mathrm{Alg}(H,Z(B^{\co H})) $.
The corresponding universal $K$-matrix for
$\alpha\in\mathrm{Alg}(H,Z(B^{\co H}))$
is given by $(\id_H\otimes F_B)(K_\alpha)$, where
\begin{equation*}
	K_\alpha=(\id_H\otimes j_\alpha)(R_{21}R),
	\quad
	j_\alpha=(\alpha\otimes\id_H)\circ\Delta^{\op},
\end{equation*}
and $F_B:B^{\co H}\otimes H\to B$ is the isomorphism
of left $H$-comodule algebras constructed from $(B,K)$
in \cref{thm:factorizable-B-structure}.
\end{theorem}

\begin{remark}\label{rem:K-classification-basepoint}
	The bijection in \cref{thm:K-matrices-general-B} is obtained after fixing a universal $K$-matrix $K$ on $B$.
	Indeed, $\mathrm{Alg}(H,Z(B^{\co H}))$ always contains the algebra map $h\mapsto\varepsilon(h)1_B$, whereas a left $H$-comodule algebra need not admit a universal $K$-matrix.
	Under this bijection, the algebra map $h\mapsto\varepsilon(h)1_B$ corresponds to $K$.
\end{remark}

We establish the following two lemmas for an arbitrary nonzero finite-dimensional algebra $A$.
We regard $A\otimes H$ as a left $H$-comodule algebra with the natural coaction given by $\delta(a\otimes h)=h_{(1)}\otimes(a\otimes h_{(2)})$.

\begin{lemma}\label{lem:colinear-algebra-maps-standard-form}
	Regard $H$ as a left $H$-comodule algebra via its coproduct.
	Let $j\in\HAlg(H,A\otimes H)$ and $\alpha\in\mathrm{Alg}(H,A)$.
	Then the assignments
	\begin{alignat}{2}
		\label{eq:alpha-from-j}
		j &\mapsto \alpha_j &&:= (\id_A\otimes\varepsilon)\circ j,\\
		\label{eq:j-alpha}
		\alpha &\mapsto j_\alpha &&:= (\alpha\otimes\id_H)\circ\Delta^{\op}
	\end{alignat}
	give mutually inverse bijections between $\HAlg(H,A\otimes H)$ and $\mathrm{Alg}(H,A)$.
\end{lemma}

\begin{proof}
	Let $j\in\HAlg(H,A\otimes H)$ and define $\alpha_j$ by \eqref{eq:alpha-from-j}.
	Clearly, $\alpha_j$ is an algebra map.
	For $h\in H$, let $j(h)=\sum_i c_i\otimes x_i$, where $c_i\in A$ and $x_i\in H$.
	Since $j$ is a left $H$-comodule map, we have
	\begin{equation*}
		(\delta \circ j) (h)
		=
		\sum_i
		x_{i(1)}\otimes c_i\otimes x_{i(2)}
		=
		h_{(1)}\otimes j(h_{(2)})
		=
		(\id_H\otimes j) (\Delta (h)).
	\end{equation*}
	Let $\tau_{H,A}:H\otimes A\to A\otimes H$ denote the flip map.
	Applying $\tau_{H,A}\circ(\id_H\otimes\id_A\otimes\varepsilon)$ to the preceding identity gives
	\begin{equation*}
		j(h)
		=
		\alpha_j(h_{(2)})\otimes h_{(1)} = j_{\alpha_j} (h).
	\end{equation*}

	Conversely, let $\alpha\in\mathrm{Alg}(H,A)$ and define $j_\alpha$ by \eqref{eq:j-alpha}.
	For all $h,k\in H$, we have
	\begin{equation*}
		j_\alpha(hk)
		=
		\alpha(h_{(2)}k_{(2)})
		\otimes
		h_{(1)}k_{(1)}
		=
		j_\alpha(h)j_\alpha(k).
	\end{equation*}
	Hence $j_\alpha$ is an algebra map.
	Moreover,
	\begin{align*}
		\delta(j_\alpha(h))
		=
		h_{(1)}
		\otimes
		(
		\alpha(h_{(3)})\otimes h_{(2)}
		)
		=
		(\id_H\otimes j_\alpha)\Delta(h),
	\end{align*}
	so $j_\alpha$ is a left $H$-comodule algebra map.
	We also have
	\begin{align*}
		\alpha_{j_{\alpha}} (h)
		=
		(\id_A
		\otimes\varepsilon)(j_\alpha(h))
		=
		\alpha(h_{(2)})\varepsilon(h_{(1)})
		=
		\alpha(h).
	\end{align*}
	Consequently, $j_{\alpha_j}=j$ and $\alpha_{j_\alpha}=\alpha$.
	Thus \eqref{eq:alpha-from-j} and \eqref{eq:j-alpha} give mutually inverse bijections.
\end{proof}

For $\alpha\in\mathrm{Alg}(H,A)$, let $j_\alpha$ be the left $H$-comodule algebra map defined by \eqref{eq:j-alpha}, and set
\begin{equation}\label{eq:K-alpha-definition}
	K_\alpha
	:=
	(\id_H\otimes j_\alpha)(R_{21}R)
	\in
	H\otimes A\otimes H.
\end{equation}

\begin{lemma}\label{lem:K-alpha-central}
	Let $\alpha\in\mathrm{Alg}(H,A)$.
	Then $K_\alpha$ is a universal $K$-matrix for $A\otimes H$ if and only if
	\begin{equation}\label{eq:alpha-central}
		\alpha(H)\subseteq Z(A).
	\end{equation}
\end{lemma}

\begin{proof}
	Suppose first that $K_\alpha$ is a universal $K$-matrix for $A\otimes H$.
	Let $j_{A\otimes H}\in\HAlg(H,A\otimes H)$ be the map associated with $(A\otimes H,K_\alpha)$ by \eqref{eq:jB}.
	Then \eqref{eq:KB-slice-QH} gives
	$j_{A\otimes H}(Q_H(f))
	=
	(f\otimes\id_{A\otimes H})(K_\alpha)$
	for all $f\in H^*$.
	We also have
	$(f\otimes\id_{A\otimes H})(K_\alpha)=j_\alpha(Q_H(f))$
	by the definition of $K_\alpha$.
	It follows that
	\begin{equation*}
		j_{A\otimes H}(Q_H(f))
		=
		j_\alpha(Q_H(f)).
	\end{equation*}
	Since $H$ is factorizable, $Q_H$ is surjective, and hence $j_{A\otimes H}=j_\alpha$.
	By \cref{lem:jB-commutes-coinvariants}, every element of $j_\alpha(H)$ commutes with every element of $A\otimes1_H=(A\otimes H)^{\co H}$.
	In particular, we have
	\begin{equation*}
		j_\alpha(h)(a\otimes1_H)
		=
		(a\otimes1_H)j_\alpha(h)
	\end{equation*}
	for all $h\in H$ and $a\in A$.
	We obtain $\alpha(h)a=a\alpha(h)$ by applying $\id_A\otimes\varepsilon$ to both sides of the preceding equality and using $j_\alpha(h)=\alpha(h_{(2)})\otimes h_{(1)}$.
	This proves \eqref{eq:alpha-central}.

	Conversely, assume \eqref{eq:alpha-central}.
	Consider the map
	\begin{equation*}
		\Gamma_\alpha:
		A\otimes H
		\to
		A\otimes H,
		\quad
		\Gamma_\alpha(a\otimes h)
		=
		a\alpha(h_{(2)})\otimes h_{(1)}.
	\end{equation*}
	Since $\alpha(H)\subseteq Z(A)$, the following identities hold for all $a,b\in A$ and $h,k\in H$:
	\begin{align*}
		\Gamma_\alpha((a\otimes h)(b\otimes k))
		&=
		ab\alpha(h_{(2)}k_{(2)})
		\otimes h_{(1)}k_{(1)} \\
		&=
		a\alpha(h_{(2)})
		b\alpha(k_{(2)})
		\otimes h_{(1)}k_{(1)}
		=
		\Gamma_\alpha(a\otimes h)
		\Gamma_\alpha(b\otimes k),
	\end{align*}
	and
	\begin{equation*}
		\delta(\Gamma_\alpha(a\otimes h))
		=
		h_{(1)}
		\otimes
		(a\alpha(h_{(3)})\otimes h_{(2)})
		=
		(\id_H\otimes\Gamma_\alpha)
		\delta(a\otimes h).
	\end{equation*}
	Thus $\Gamma_\alpha$ is a left $H$-comodule algebra map.
	Moreover, for all $a\in A$ and $h\in H$, we have
	\begin{equation*}
		\Gamma_\alpha
		(a\alpha(S^{-1}(h_{(2)}))\otimes h_{(1)})
		=
		a\alpha(S^{-1}(h_{(3)})h_{(2)})\otimes h_{(1)}
		=
		a\otimes h,
	\end{equation*}
	which shows that $\Gamma_\alpha$ is surjective.
	Since $(R_{21}R)_{13}$ is a universal $K$-matrix for $A\otimes H$, \cref{prop:K-matrix-quotient} implies that
	\begin{align*}
		(\id_H\otimes\Gamma_\alpha)((R_{21}R)_{13})
		=
		(\id_H\otimes j_\alpha)(R_{21}R)
		=
		K_\alpha
	\end{align*}
	is a universal $K$-matrix for $A\otimes H$.
\end{proof}

\begin{proposition}\label{prop:K-matrices-tensor-H}
	The map $\alpha\mapsto K_\alpha$, with $K_\alpha$ as in \eqref{eq:K-alpha-definition}, defines a bijection from $\mathrm{Alg}(H,Z(A))$ to the set of universal $K$-matrices on $A\otimes H$.
\end{proposition}

\begin{proof}
	By the initiality of $(H,R_{21}R)$, every universal $K$-matrix on $A\otimes H$ is of the form $(\id_H\otimes j)(R_{21}R)$ for a unique $j\in\HAlg(H,A\otimes H)$.
	By \cref{lem:colinear-algebra-maps-standard-form}, every such $j$ equals $j_\alpha$ for a unique $\alpha\in\mathrm{Alg}(H,A)$.
	The assertion now follows from \cref{lem:K-alpha-central}.
\end{proof}

The correspondence in \cref{prop:K-matrices-tensor-H} may be compared with the characterization of algebra maps from the reflection equation dual in \cite[Proposition~10]{DKM03} and the description of braided module structures in terms of algebra maps from the reflection equation algebra in \cite[Theorem~4.9 and Corollary~4.11]{BZBJ18}.
Here we work with universal $K$-matrices on left $H$-comodule algebras and obtain an explicit classification for $A\otimes H$ in terms of $\mathrm{Alg}(H,Z(A))$.
See also \cite[Proposition~3.13]{Saf21} for a related characterization of quantum moment maps by centrality in a braided category.

\begin{proof}[Proof of \cref{thm:K-matrices-general-B}]
	Let $F_B:B^{\co H}\otimes H\to B$ be the isomorphism of left $H$-comodule algebras given in \cref{thm:factorizable-B-structure}.
	The isomorphism $F_B$ induces a bijection between the universal $K$-matrices on $B^{\co H}\otimes H$ and those on $B$ via $\id_H\otimes F_B$.
	Applying \cref{prop:K-matrices-tensor-H} with $A=B^{\co H}$ gives the desired bijection sending $\alpha\in\mathrm{Alg}(H,Z(B^{\co H}))$ to $(\id_H\otimes F_B)(K_\alpha)$.
\end{proof}

\subsection{Bijectivity of \texorpdfstring{$\theta_B$}{theta_B} and factorizability}
\label{subsec:factorizability-over-factorizable-H}

When $(H,R)$ is factorizable, we characterize the bijectivity of $\theta_B:H^*\to E(H,B)$ without assuming the exactness or $H$-indecomposability of $B$.
We then compare this condition with the factorizability of $(B,K)$.
Throughout this subsection, we retain the notation introduced above and in \cref{subsec:factorizability-quasitriangular-comodule-algebras}.

We first establish a criterion for the bijectivity of $\theta_B$.

\begin{theorem}
\label{thm:theta-B-bijective-center}
Let $(H,R)$ be a finite-dimensional factorizable Hopf algebra and $(B,K)$ a finite-dimensional quasitriangular left $H$-comodule algebra.
Then the following conditions are equivalent:
\begin{enumerate}
    \item[\rm{(i)}] The map $\theta_B:H^*\to E(H,B)$ is bijective.
    \item[\rm{(ii)}] $\dim E(H,B)=\dim H$.
    \item[\rm{(iii)}] $Z(B^{\co H})=\bfk$.
\end{enumerate}
\end{theorem}

We now prove that $\theta_B$ is injective.
The injectivity of $\theta_{R_H(A)}$ is known without the factorizability assumption on $(H,R)$ \cite[Proof of Proposition~6.24]{WY25}.

\begin{proposition}
\label{prop:theta-injective-factorizable-H}
Suppose that $(H,R)$ is factorizable.
Then the linear map $\theta_B:H^*\to E(H,B)$ defined by \eqref{eq:theta-B} is injective for every finite-dimensional quasitriangular left $H$-comodule algebra $(B,K)$.
\end{proposition}

\begin{proof}
Let $f\in H^*$.
Setting $h=1_H$ in \eqref{eq:theta-B} gives
\begin{equation*}
    \theta_B(f)(1_H)
    =
    (f\otimes\id_B)(K).
\end{equation*}
Together with \eqref{eq:KB-slice-QH}, this yields $\theta_B(f)(1_H)=j_B(Q_H(f))$.
Suppose that $\theta_B(f)=0$.
Then $j_B(Q_H(f))=0$.
Let $F_B:B^{\co H}\otimes H\to B$ be the isomorphism of quasitriangular left $H$-comodule algebras from \cref{thm:factorizable-B-structure}.
The equality $j_B(h)=F_B(1_B\otimes h)$ holds for all $h\in H$.
The map $H\to B^{\co H}\otimes H$, $h\mapsto1_B\otimes h$, is injective, and so is $j_B$ because $F_B$ is an isomorphism of left $H$-comodule algebras.
Since $(H,R)$ is factorizable, $Q_H:H^*\to H$ is a vector space isomorphism.
It follows that $f=0$.
\end{proof}

\begin{lemma}
\label{lem:E-general-factorizable}
There is an isomorphism of vector spaces $E(H,B)\cong Z(B^{\co H})\otimes H$.
\end{lemma}

\begin{proof}
By \cref{thm:factorizable-B-structure}, composition with $F_B$ gives a vector space isomorphism from $E(H,B^{\co H}\otimes H)$ to $E(H,B)$.
It therefore suffices to construct an isomorphism of vector spaces from $E(H,B^{\co H}\otimes H)$ to $Z(B^{\co H})\otimes H$.

Let $\xi\in E(H,B^{\co H}\otimes H)$.
Applying \eqref{eq:E-H-B} to $b=c\otimes1_H$ shows that $\xi(h)$ commutes with $c\otimes1_H$ for all $c\in B^{\co H}$ and $h\in H$.
Thus $\xi(h)\in Z(B^{\co H})\otimes H$ for all $h\in H$.
In particular, evaluation at $1_H$ defines a linear map
\begin{equation*}
    \mathrm{ev}_{1_H}:
    E(H,B^{\co H}\otimes H)
    \to
    Z(B^{\co H})\otimes H,
    \quad
    \xi\mapsto\xi(1_H).
\end{equation*}
We show that this map is bijective.

For $g\in H$, the defining relation \eqref{eq:E-H-B} with $b=1_B\otimes g$ and $h=1_H$ gives
\begin{equation*}
    \xi(g_{(1)})(1_B\otimes g_{(2)})
    =
    (1_B\otimes g)\xi(1_H).
\end{equation*}
The antipode identities then give
\begin{equation*}
    \xi(g)
    =
    \xi(g_{(1)})(1_B\otimes g_{(2)}S(g_{(3)}))
    =
    (1_B\otimes g_{(1)})\xi(1_H)(1_B\otimes S(g_{(2)})).
\end{equation*}
If $\xi(1_H)=0$, this equality implies that $\xi=0$.
Hence $\mathrm{ev}_{1_H}$ is injective.
To prove surjectivity, let $y\in Z(B^{\co H})\otimes H$ and define the linear map $\xi_y:H\to B^{\co H}\otimes H$ by
\begin{equation*}
    \xi_y(g)
    =
    (1_B\otimes g_{(1)})y(1_B\otimes S(g_{(2)})).
\end{equation*}
The definition gives $\xi_y(1_H)=y$.
It remains to verify that $\xi_y\in E(H,B^{\co H}\otimes H)$.
The following identity holds for all $c\in B^{\co H}$ and $g,h\in H$:
\begin{align*}
    \xi_y(g_{(1)}h)(c\otimes g_{(2)})
    &=
    (1_B\otimes g_{(1)}h_{(1)})y(1_B\otimes S(h_{(2)})S(g_{(2)}))(c\otimes g_{(3)})
    \\
    &=
    (1_B\otimes gh_{(1)})y(c\otimes S(h_{(2)}))
    \\
    &=
    (c\otimes gh_{(1)})y(1_B\otimes S(h_{(2)}))
    =
    (c\otimes g)\xi_y(h).
\end{align*}
The second equality follows from the antipode identity $S(g_{(2)})g_{(3)}=\varepsilon(g_{(2)})1_H$, and the third uses $y(c\otimes1_H)=(c\otimes1_H)y$.
Thus $\xi_y\in E(H,B^{\co H}\otimes H)$, and $\mathrm{ev}_{1_H}(\xi_y)=y$ proves surjectivity.
We conclude that $\mathrm{ev}_{1_H}$ is an isomorphism of vector spaces.
\end{proof}

\begin{proof}[Proof of \cref{thm:theta-B-bijective-center}]
The map $\theta_B$ is injective by \cref{prop:theta-injective-factorizable-H}, so it is bijective if and only if $\dim E(H,B)=\dim H^*=\dim H$.
This proves the equivalence of \textup{(i)} and \textup{(ii)}.
The isomorphism in \cref{lem:E-general-factorizable} gives
\begin{equation*}
    \dim E(H,B)
    =
    \dim Z(B^{\co H})\cdot\dim H.
\end{equation*}
As $1_B\in Z(B^{\co H})$, the condition $\dim Z(B^{\co H})=1$ is equivalent to $Z(B^{\co H})=\bfk$.
This establishes the equivalence of \textup{(ii)} and \textup{(iii)}.
\end{proof}

We next characterize $H$-simplicity of $B$ and use it to obtain a criterion for factorizability of $(B,K)$.

\begin{proposition}
\label{prop:H-simple-coinvariants}
Let $(H,R)$ be a factorizable Hopf algebra and $(B,K)$ a quasitriangular $H$-comodule algebra.
Then $B$ is $H$-simple if and only if $B^{\co H}$ is simple as a $\bfk$-algebra.
\end{proposition}

\begin{proof}
Let $F_B:B^{\co H}\otimes H\to B$ be the isomorphism of left $H$-comodule algebras from \cref{thm:factorizable-B-structure}.

Let $I\subseteq B^{\co H}$ be an ideal of $B^{\co H}$.
The subspace $I\otimes H$ is an ideal of $B^{\co H}\otimes H$ and is stable under the left $H$-coaction \eqref{eq:coinvariants-tensor-H-coaction}.
Its image $F_B(I\otimes H)$ is therefore an $H$-costable ideal of $B$.

Now let $J\subseteq B$ be an $H$-costable ideal of $B$.
The right $H$-action \eqref{eq:right-action-on-B} preserves $J$ because $J$ is an ideal of $B$.
Moreover, $J^{\co H}=J\cap B^{\co H}$.
The Fundamental Theorem of Hopf Modules \cite[Theorem~1.9.4]{Mon93} gives an isomorphism of vector spaces $(J\cap B^{\co H})\otimes H\to J$, $x\otimes h\mapsto xj_B(h)$.
This is the restriction of $F_B$ to $(J\cap B^{\co H})\otimes H$, so $  F_B^{-1}(J)
=
(J\cap B^{\co H})\otimes H. $
The intersection $J\cap B^{\co H}$ is an ideal of $B^{\co H}$.
Therefore the assignments
\begin{equation*}
    I\mapsto F_B(I\otimes H),
    \quad
    J\mapsto J\cap B^{\co H}
\end{equation*}
give mutually inverse bijections between the ideals of $B^{\co H}$ and the $H$-costable ideals of $B$.
Hence, $B$ is $H$-simple if and only if $B^{\co H}$ is simple.
\end{proof}

The preceding results yield the following criterion.

\begin{proposition}
\label{prop:factorizable-B-coinvariants}
With the notation above, the following conditions are equivalent:
\begin{enumerate}[label=\rm{(\roman*)}]
    \item $(B,K)$ is a factorizable left $H$-comodule algebra.
    \item $B$ is $H$-simple.
    \item $B^{\co H}$ is simple as a $\bfk$-algebra.
\end{enumerate}
\end{proposition}

\begin{proof}
The equivalence of \textup{(i)} and \textup{(ii)} follows from \cref{prop:H-simple-exact-indecomposable,prop:theta-injective-factorizable-H} and \cref{lem:factorizable-comodule-injective}.
The equivalence of \textup{(ii)} and \textup{(iii)} is precisely \cref{prop:H-simple-coinvariants}.
\end{proof}

\begin{remark}
\label{rem:theta-bijective-factorizable}
The map $\theta_B:H^*\to E(H,B)$ is bijective if and only if $Z(B^{\co H})=\bfk$ by \cref{thm:theta-B-bijective-center}, but $B$ is $H$-simple if and only if $B^{\co H}$ is simple by \cref{prop:H-simple-coinvariants}.
In particular, $H$-simplicity of $B$ implies bijectivity of $\theta_B$, but the converse does not hold in general (see \cref{ex:exm-reflective-theta-bijective-not-H-simple}).
\end{remark}

\section{Reflective algebras over factorizable Hopf algebras}
\label{sec:factorizable-reflective-case}

Throughout this section, let $(H,R)$ be a finite-dimensional quasitriangular Hopf algebra and $A$ a nonzero finite-dimensional left $H$-comodule algebra.
We apply the general results on quasitriangular comodule algebras from the preceding section to the reflective algebra $(R_H(A),K_{\mathrm{ref}}(A))$, using the results of \cref{sec:drinfeld-quotient}.

\subsection{Structure and universal $K$-matrices of reflective algebras}
\label{subsec:factorizable-reflective-algebras}

In this subsection, assume that $(H,R)$ is factorizable.
We use \cref{sec:drinfeld-quotient} to describe explicitly the results obtained by applying \cref{thm:factorizable-B-structure} to $B=R_H(A)$.

Since $H$ is factorizable, we have $\Ker Q_H=0$, $\Img Q_H=H$, and $T=\bfk$.
Consequently, the isomorphism in \cref{thm:cotensor-quasitriangular} becomes the following isomorphism of left $H$-comodule algebras:
\begin{equation}\label{eq:cotensor-factorizable}
R_H(A)
\xrightarrow{\sim}
H\otimes A,
\quad
a\xi
\mapsto
a_{(-1)}Q_H(\xi)\otimes a_{(0)}.
\end{equation}
Here $H\otimes A$ is equipped with the left $H$-coaction $\delta_{H\otimes A}(h\otimes a)=h_{(1)}\otimes(h_{(2)}\otimes a)$.
We also equip $A\otimes H$ with the left $H$-coaction $\delta_{A\otimes H}(a\otimes h)=h_{(1)}\otimes(a\otimes h_{(2)})$.
The flip map
\begin{equation*}
\tau_{H,A}:H\otimes A\to A\otimes H,
\quad
h\otimes a\mapsto a\otimes h
\end{equation*}
is an isomorphism of left $H$-comodule algebras.
Let $A_{\mathrm{triv}}$ denote the algebra $A$ equipped with the trivial left $H$-coaction $a\mapsto 1_H\otimes a$.

\begin{theorem}
\label{thm:factorizable-reflective-structure}
With the notation above, the map
\begin{equation*}
\Psi_A:
R_H(A)
\to
A_{\mathrm{triv}}\otimes H,
\quad
\Psi_A(a\xi)
=
a_{(0)}\otimes a_{(-1)}Q_H(\xi)
\end{equation*}
induces an isomorphism of quasitriangular left $H$-comodule algebras
\begin{equation*}
(R_H(A),K_{\mathrm{ref}}(A))
\cong
(A_{\mathrm{triv}}\otimes H,(R_{21}R)_{13}).
\end{equation*}
\end{theorem}

\begin{proof}
The map $\Psi_A$ is the composition of the isomorphism \eqref{eq:cotensor-factorizable} and the flip map $\tau_{H,A}$, so it is an isomorphism of left $H$-comodule algebras.
The canonical universal $K$-matrix $\overline{K}_{\mathrm{ref}}(A)$ on $R_H(A)/A\,\Ker Q_H$ corresponds to $(R_{21}R)_{12}$ on $H\square_T A$ by \cref{thm:cotensor-quasitriangular}.
Since $\Ker Q_H=0$ and $T=\bfk$, the quotient $R_H(A)/A\,\Ker Q_H$ is $R_H(A)$ and $\overline{K}_{\mathrm{ref}}(A)=K_{\mathrm{ref}}(A)$.
Thus $K_{\mathrm{ref}}(A)$ corresponds to $(R_{21}R)_{12}$ on $H\otimes A$.
Applying $\id_H\otimes\tau_{H,A}$ sends $(R_{21}R)_{12}$ to $(R_{21}R)_{13}$.
Therefore, $\Psi_A$ is an isomorphism of quasitriangular left $H$-comodule algebras.
\end{proof}

\begin{corollary}
\label{cor:RH-A-coinvariants}
The restriction of $\Psi_A$ to $R_H(A)^{\co H}$ yields an isomorphism of algebras $R_H(A)^{\co H}\cong A$.
\end{corollary}

\begin{proof}
The map $\Psi_A$ is an isomorphism of left $H$-comodule algebras, so it restricts to an isomorphism of algebras between their coinvariant subalgebras.
In particular, $R_H(A)^{\co H}\cong (A_{\mathrm{triv}}\otimes H)^{\co H}=A\otimes1_H\cong A$ as algebras.
\end{proof}
It follows from \eqref{eq:coinvariants-tensor-H-coaction} that the left $H$-coaction on $A_{\mathrm{triv}}\otimes H$ is given by
\begin{equation}\label{eq:coaction-of-Atriv-otimes-H}
	\delta(a\otimes h)=h_{(1)}\otimes(a\otimes h_{(2)})\quad(a\in A,\ h\in H).
\end{equation}

\begin{theorem}
\label{thm:K-matrices-reflective}
We retain the notation above.
The universal $K$-matrices on $R_H(A)$ are in bijection with the algebra maps $\alpha\in\mathrm{Alg}(H,Z(A))$.
The corresponding universal $K$-matrix for $\alpha\in\mathrm{Alg}(H,Z(A))$ is given by $ (\id_H\otimes\Psi_A^{-1})
(\id_H\otimes j_\alpha)(R_{21}R) $,
where $j_\alpha:H\to A_{\mathrm{triv}}\otimes H$ is the algebra map defined by
$j_\alpha=(\alpha\otimes\id_H)\circ\Delta^{\op}$.
\end{theorem}

\begin{proof}
We may apply \cref{thm:K-matrices-general-B} with $B=R_H(A)$ by \cref{thm:factorizable-reflective-structure,cor:RH-A-coinvariants}.
Since $R_H(A)^{\co H}\cong A$ as algebras, their centers are isomorphic.
The stated bijection follows from \cref{thm:K-matrices-general-B}.
The formula follows by applying the same theorem to $(A_{\mathrm{triv}}\otimes H,(R_{21}R)_{13})$ and then applying $\id_H\otimes\Psi_A^{-1}$ to the corresponding universal $K$-matrix.
\end{proof}

If $A$ is commutative, \cref{thm:K-matrices-reflective} gives a bijection between the universal $K$-matrices on $R_H(A)$ and $\mathrm{Alg}(H,A)$.
Taking $A=\bfk$ gives the following corollary.
\begin{corollary}
\label{cor:K-matrices-H}
Let $(H,R)$ be a factorizable Hopf algebra.
The universal $K$-matrices on the quasitriangular left $H$-comodule algebra $(H,R_{21}R)$ are in bijection with $\mathrm{Alg}(H,\bfk)$.
More explicitly, for $\chi\in\mathrm{Alg}(H,\bfk)$, define the map $j_\chi:H\to H$ by $j_\chi(h)=\chi(h_{(2)})h_{(1)}$.
The corresponding universal $K$-matrix is given by
\begin{equation*}
K_\chi
=
(\id_H\otimes j_\chi)(R_{21}R).
\end{equation*}
In particular, $K_\varepsilon=R_{21}R$.
\end{corollary}

\begin{proof}
By \cref{cor:factorizable-H-drinfeld-map}, $(H,R_{21}R)$ and $(R_H(\bfk),K_{\mathrm{ref}}(\bfk))$ are isomorphic as quasitriangular left $H$-comodule algebras.
The bijection follows from \cref{thm:K-matrices-reflective}, and its explicit description is given by \cref{lem:colinear-algebra-maps-standard-form}.
\end{proof}

\subsection{Factorizability of reflective algebras}
\label{subsec:theta-reflective-bijective}

We first assume that $(H,R)$ is factorizable and apply the results of \cref{subsec:factorizability-over-factorizable-H} to $(R_H(A),K_{\mathrm{ref}}(A))$ to characterize the bijectivity of $\theta_{R_H(A)}$.
We then use \cref{subsec:costable-ideals-simplicity} to obtain necessary and sufficient conditions on $H$ and $A$ for $R_H(A)$ to be $H$-simple without assuming that $(H,R)$ is factorizable.
Combining this characterization with the results of \cite{WY25} gives a criterion for the factorizability of reflective algebras.

\begin{theorem}
\label{thm:theta-reflective-bijective-center}
Suppose that $(H,R)$ is factorizable.
Then the following conditions are equivalent.
\begin{enumerate}
\item[\rm{(i)}] The map $\theta_{R_H(A)}:H^*\to E(H,R_H(A))$ is bijective.
\item[\rm{(ii)}] $\dim E(H,R_H(A))=\dim H$.
\item[\rm{(iii)}] $Z(A)=\bfk$.
\end{enumerate}
\end{theorem}

\begin{proof}
The result follows by applying \cref{thm:theta-B-bijective-center} with $B=R_H(A)$ and using \cref{cor:RH-A-coinvariants}.
\end{proof}

We now determine when $(R_H(A),K_{\mathrm{ref}}(A))$ is factorizable without assuming that $(H,R)$ is factorizable.

\begin{proposition}
\label{prop:factorizable-H-RH-simple}
Let $(H,R)$ be a quasitriangular Hopf algebra, not necessarily factorizable, and let $A$ be a nonzero left $H$-comodule algebra.
The following conditions are equivalent.
\begin{enumerate}
\item[\rm{(i)}] The reflective algebra $R_H(A)$ is $H$-simple.
\item[\rm{(ii)}] The Hopf algebra $(H,R)$ is factorizable and $A$ is simple as an algebra.
\end{enumerate}
\end{proposition}

\begin{proof}
Recall that $(H,R)$ is factorizable if and only if $\Ker Q_H=0$.

Suppose that \textup{(i)} holds.
The subspace $A\,\Ker Q_H\subseteq R_H(A)$ is an $H$-costable ideal of $R_H(A)$ by \cref{prop:canonical-ideal}.
We also have the following isomorphisms of vector spaces:
\begin{equation*}
A\,\Ker Q_H
\cong
A\otimes\Ker Q_H,
\quad
R_H(A)/A\,\Ker Q_H
\cong
A\otimes\Img  Q_H.
\end{equation*}
Since $A\neq0$ and $Q_H(\varepsilon)=1_H$, the second isomorphism shows that the quotient is nonzero.
Thus $A\,\Ker Q_H$ is a proper ideal of $R_H(A)$.
The $H$-simplicity of $R_H(A)$ therefore forces $A\,\Ker Q_H=0$.
Since $A\neq0$, the first vector space isomorphism implies $\Ker Q_H=0$, and hence $(H,R)$ is factorizable.
Moreover, $R_H(A)/A\,\Ker Q_H=R_H(A)$ is $H$-simple, so $A$ is $T$-simple by \cref{cor:quotient-H-simple}.
Now $T=\bfk$, and the induced left $T$-coaction on $A$ is trivial.
Thus $T$-simplicity of $A$ is equivalent to simplicity of $A$ as an algebra.

Now suppose that \textup{(ii)} holds.
We then have $\Ker Q_H=0$ and $T=\bfk$.
Since $A$ is simple as an algebra, it is $T$-simple, and \cref{cor:quotient-H-simple} implies that $R_H(A)$ is $H$-simple.
\end{proof}

We obtain the following characterization of factorizability.

\begin{theorem}
\label{thm:reflective-factorizability-criterion}
The following conditions are equivalent.
\begin{enumerate}
\item[\rm{(i)}] The quasitriangular left $H$-comodule algebra $(R_H(A),K_{\mathrm{ref}}(A))$ is factorizable.
\item[\rm{(ii)}] The quasitriangular Hopf algebra $(H,R)$ is factorizable and $A$ is simple as a $\bfk$-algebra.
\end{enumerate}
In particular, if $(H,R)$ is factorizable, then $(R_H(A),K_{\mathrm{ref}}(A))$ is factorizable if and only if $A$ is simple as an algebra.
\end{theorem}

\begin{proof}
The result follows from \cref{cor:reflective-factorizable-H-simple,prop:factorizable-H-RH-simple}.
\end{proof}

\begin{remark}
\label{rem:comparison-WY-reflective}
Suppose that $(H,R)$ is factorizable.
By \cref{thm:theta-reflective-bijective-center,prop:factorizable-H-RH-simple}, the map $\theta_{R_H(A)}$ is bijective if and only if $Z(A)=\bfk$, but $R_H(A)$ is $H$-simple if and only if $A$ is simple as an algebra.
Thus the $H$-simplicity of $R_H(A)$ implies the bijectivity of $\theta_{R_H(A)}$, but the converse does not hold in general (see \cref{ex:exm-reflective-theta-bijective-not-H-simple}).
\end{remark}

\section{Examples of reflective algebras over $u_q(\mathfrak{sl}_2)$}\label{sec:example-uqsl2}

In this section, we consider the small quantum group $H=u_q(\mathfrak{sl}_2)$ associated with a root of unity $q$ of odd order as an example of our results.
We first identify $R_H(\bfk)$ with $H$ as quasitriangular left $H$-comodule algebras.
We then recall the family $A_2(r;\zeta)$ of left $H$-comodule algebras from \cite[Section~5.6]{NSS25} and investigate the structure and factorizability of their reflective algebras $R_H(A_2(r;\zeta))$.
Throughout this section, we further assume that $\bfk$ has characteristic zero.

\subsection{The factorizable Hopf algebra $u_q (\mathfrak{sl}_2)$}
\label{subsec:Hopf-algebra-u_q}

Let $N > 1$ be an odd integer, and let $q \in \bfk$ be a primitive $N$-th root of unity.
We define $u_q(\mathfrak{sl}_2)$ to be the algebra generated by $K$, $E$ and $F$ subject to the relations
\begin{gather*}
	K^N = 1,
	\ K E = q^2 E K,
	\ K F = q^{-2} F K,
	\ E^N = F^N = 0,
	\ E F - F E = \frac{K - K^{-1}}{q - q^{-1}}.
\end{gather*}
The algebra $u_q(\mathfrak{sl}_2)$ is a Hopf algebra with the comultiplication $\Delta$, the counit $\varepsilon$ and the antipode $S$ determined by
\begin{gather*}
	\Delta(K) = K \otimes K,
	\ \Delta(E) = E \otimes K + 1 \otimes E,
	\ \Delta(F) = F \otimes 1 + K^{-1} \otimes F, \\
	\ \varepsilon(K) = 1,
	\ \varepsilon(E) = \varepsilon(F) = 0,
	\ S(K) = K^{-1},
	\ S(E) = -E K^{-1},
	\ S(F) = -K F.
\end{gather*}
In what follows, we put $H=u_q(\mathfrak{sl}_2)$.

For $ q \in \bfk^{\times} $ and $m \in \mathbb{Z}_{\ge 0}$, we define
\begin{equation*}
	[m]_q= 1 + q + \cdots + q^{m-1}
	\quad \text{and} \quad
	[m]_q!=[1]_q [2]_q \cdots [m]_q
\end{equation*}
with the conventions $[0]_q = 0$ and $[0]_q! = 1$.
By \cite[Theorem~IX.7.1]{Kas95}, the following element
$R\in H\otimes H$ is a quantum $R$-matrix for $H$:
\begin{equation}\label{eq:R-uqsl2}
	R=\frac{1}{N}\sum_{m,i,j=0}^{N-1}
	\frac{(q-q^{-1})^m}{[m]_{q^2}!}
	q^{m(m-1)+2m(i-j)-2ij}
	E^mK^i\otimes F^mK^j.
\end{equation}
Moreover, $(H,R)$ is factorizable
\cite[Theorem~5.9 and Example~5.11]{LO17}.
Consequently, the Drinfeld map $Q_H:R_H(\bfk)\to H$ gives an isomorphism of quasitriangular left $H$-comodule algebras $(R_H(\bfk),K_{\mathrm{ref}}(\bfk))\cong(H,R_{21}R)$ by \cref{prop:drinfeld-map-QT-map}.

We note the following observation.
\begin{proposition}\label{prop:exm-central-character}
	Let $A$ be a nonzero $\bfk$-algebra.
	Then every algebra map $\alpha:H\to Z(A)$ is given by $\alpha(h)=\varepsilon(h)1_A$ for all $h\in H$.
\end{proposition}

\begin{proof}
	Every algebra map from $H$ to a commutative algebra factors through $H_{\mathrm{ab}}$, the abelianization of $H$.
	The relations $KE=q^2EK$ and $KF=q^{-2}FK$ imply that $E=F=0$ in $H_{\mathrm{ab}}$, while the remaining relations give $K^2=K^N=1$, and hence $K=1$ since $N$ is odd.
	Thus $H_{\mathrm{ab}}\cong\bfk$, which proves the claim.
\end{proof}

\begin{corollary}
	\label{cor:exm-unique-K-H}
	The left $H$-comodule algebra $H$ with coaction $\Delta$ has a unique universal $K$-matrix, namely $R_{21}R$.
\end{corollary}

\begin{proof}
	We have $\mathrm{Alg} (H,\bfk)=\{\varepsilon\}$ by \cref{prop:exm-central-character}.
	It follows from \cref{cor:K-matrices-H} that $R_{21}R$ is the unique universal $K$-matrix on $H$.
\end{proof}

\subsection{\texorpdfstring{$u_q(\mathfrak{sl}_2)$}{u_q(sl_2)}-comodule algebras \texorpdfstring{$A_2(r;\zeta)$}{A_2(r;zeta)}}
\label{subsec:exm-comod-alg-uq}

We recall the family $A_2(r;\zeta)$ of left $H$-comodule algebras appearing in the classification of indecomposable exact module categories over $H\text{-}\mathrm{mod}$ up to $H$-equivariant Morita equivalence \cite[Section~5.6]{NSS25}.

Let $r>1$ be a divisor of $N$, and let $\zeta\in\bfk$.
The algebra $A_2(r;\zeta)$ is generated by $G$ and $Y$ subject to
\begin{equation*}
    G^r=1,\quad Y^N=\zeta,\quad GY=q^{-2N/r}YG.
\end{equation*}
The left $H$-coaction $\delta_{A_2}$ on $A_2(r;\zeta)$ is determined by
\begin{equation*}
    \delta_{A_2}(G)=K^{N/r}\otimes G,\quad
    \delta_{A_2}(Y)=F\otimes1+K^{-1}\otimes Y.
\end{equation*}
The dimension of $A_2(r;\zeta)$ is $Nr$.

Suppose that $A_2(r;\zeta)$ has a universal $K$-matrix.
Then $\dim H=N^3$ divides $\dim A_2(r;\zeta)=Nr$ by \cref{cor:K-matrix-existence}.
However, this is impossible since $r\leq N$.
Therefore, $A_2(r;\zeta)$ has no universal $K$-matrix.

\subsection{Explicit description of \texorpdfstring{$R_H(A_2(r;\zeta))$}{R_H(A_2(r;zeta))}}
\label{subsec:exm-structure-reflective-algebras}

In this subsection, we explicitly give generators, relations, and the left $H$-coaction for the reflective algebra associated with $A_2(r;\zeta)$ introduced in \cref{subsec:exm-comod-alg-uq}.

Let $A=A_2(r;\zeta)$.
Since $(H,R)$ is factorizable, \cref{thm:factorizable-reflective-structure} gives an isomorphism of quasitriangular left $H$-comodule algebras
\begin{align*}
    \Psi_A:(R_H(A),K_{\mathrm{ref}}(A)) &\xrightarrow{\sim} (A_{\mathrm{triv}}\otimes H,(R_{21}R)_{13}),\\
    a\xi' &\mapsto a_{(0)}\otimes a_{(-1)}Q_H(\xi') \quad (a\in A,\ \xi'\in R_H(\bfk)).
\end{align*}
We identify $A$ and $R_H(\bfk)$ with their images in $R_H(A)$ under the canonical embeddings.
Since $Q_H:R_H(\bfk)\to H$ is an isomorphism of left $H$-comodule algebras, we define $k,e,f\in R_H(\bfk)$ by $k=Q_H^{-1}(K)$, $e=Q_H^{-1}(E)$ and $f=Q_H^{-1}(F)$.
By the definition of $\Psi_A$, we have
\begin{equation*}
	\Psi_A(k)=1\otimes K,\quad
	\Psi_A(e)=1\otimes E,\quad
	\Psi_A(f)=1\otimes F.
\end{equation*}
Next, we define the corresponding elements $g,y\in R_H(A)$ by
$ g=Gk^{-N/r},\  y=(Y-f)k $.
The coactions on the generators and the definition of $\Psi_A$ give
\begin{align*}
	\Psi_A(g)
	&=\Psi_A(Gk^{-N/r})
	=(G\otimes K^{N/r})(1\otimes K^{-N/r})
	=G\otimes1,\\
	\Psi_A(y)
	&=\Psi_A((Y-f)k)
	=(Y\otimes K^{-1}+1\otimes F-1\otimes F)(1\otimes K)
	=Y\otimes1.
\end{align*}
We obtain the following presentation of $R_H(A)$.

\begin{proposition}\label{prop:exm-reflective-presentation}
	The algebra $R_H(A)$ is generated by $k,e,f,g,y$ subject to the defining relations
	\begin{gather*}
		k^N=1,\quad e^N=f^N=0,\quad
		ke=q^2ek,\quad kf=q^{-2}fk,\quad
		ef-fe=\frac{k-k^{-1}}{q-q^{-1}},\\
		g^r=1,\quad y^N=\zeta,\quad
		gy=q^{-2N/r}yg,\\
		kg=gk,\quad ky=yk,\quad
		eg=ge,\quad ey=ye,\quad
		fg=gf,\quad fy=yf.
	\end{gather*}
	The left $H$-coaction $\delta_{\mathrm{ref}}$ is determined by
	\begin{gather*}
		\delta_{\mathrm{ref}}(k)=K\otimes k,\quad
		\delta_{\mathrm{ref}}(e)=E\otimes k+1\otimes e,\quad
		\delta_{\mathrm{ref}}(f)=F\otimes1+K^{-1}\otimes f,\\
		\delta_{\mathrm{ref}}(g)=1\otimes g,\quad
		\delta_{\mathrm{ref}}(y)=1\otimes y.
	\end{gather*}
\end{proposition}

\begin{proof}
The isomorphism $\Psi_A:R_H(A)\to A_{\mathrm{triv}}\otimes H$ sends $k,e,f$ to $1\otimes K,1\otimes E,1\otimes F$, respectively.
The map $\Psi_A$ also sends $g,y$ to $G\otimes1,Y\otimes1$, respectively.
The stated relations are therefore precisely those of the tensor product algebra $A_{\mathrm{triv}}\otimes H$.
The coaction formulas follow from \eqref{eq:coaction-of-Atriv-otimes-H}.
\end{proof}

\begin{corollary}\label{cor:exm-reflective-unique-K}
For every $A=A_2(r;\zeta)$, the reflective algebra $R_H(A)$ has a unique universal $K$-matrix, namely $K_{\mathrm{ref}}(A)=(\id_H\otimes Q_H^{-1})(R_{21}R)$.
Here $R$ is the universal quantum $R$-matrix given in \eqref{eq:R-uqsl2}.
We regard $Q_H^{-1}(h)\in R_H(\bfk)$ as an element of $R_H(A)$ via the canonical embedding for every $h\in H$.
\end{corollary}

\begin{proof}
The universal $K$-matrices on $R_H(A)$ are in bijection with $\mathrm{Alg} (H,Z(A))$ by \cref{thm:K-matrices-reflective}.
By \cref{prop:exm-central-character}, we have
$\mathrm{Alg}(H,Z(A))=\{\varepsilon\}$.
Thus $K_{\mathrm{ref}}(A)$ is the unique universal $K$-matrix on $R_H(A)$.
\end{proof}

\subsection{Factorizability of \texorpdfstring{$R_H(A_2(r;\zeta))$}{R_H(A_2(r;zeta))}}
\label{subsec:exm-H-simplicity-factorizability}

In this subsection, we first summarize the factorizability of $H$ and $R_H(H)$.
We then determine when $A_2(r;\zeta)$ is simple as an algebra and characterize the $H$-simplicity and factorizability of $R_H(A_2(r;\zeta))$.
Finally, using $A_2(r;\zeta)$, we give an example in which $\theta_{R_H(A_2(r;\zeta))}$ is bijective but $R_H(A_2(r;\zeta))$ is not $H$-simple, as discussed in \cref{rem:comparison-WY-reflective}.

We first consider the $H$-simplicity and factorizability of $H$ itself.
By \cref{cor:factorizable-H-drinfeld-map}, $(R_H(\bfk),K_{\mathrm{ref}}(\bfk))$ and $(H,R_{21}R)$ are isomorphic as quasitriangular left $H$-comodule algebras.
Since $H^{\co H}=\bfk$, the algebra $H$ is $H$-simple by \cref{prop:H-simple-coinvariants}.
It follows from \cref{prop:factorizable-B-coinvariants} that $(R_H(\bfk),K_{\mathrm{ref}}(\bfk))\cong(H,R_{21}R)$ is a factorizable quasitriangular left $H$-comodule algebra.
On the other hand, $H$ is not simple as an algebra, so $(R_H(H),K_{\mathrm{ref}}(H))$ is not factorizable.
We summarize these observations in the following proposition.

\begin{proposition}\label{prop:factorizability-H-and-RH(H)}
The quasitriangular left $H$-comodule algebra $(H,R_{21}R)$ is factorizable, but the reflective algebra $(R_H(H),K_{\mathrm{ref}}(H))$ is not factorizable.
\end{proposition}

We next determine when $A_2(r;\zeta)$ is simple as an algebra.

\begin{proposition}\label{prop:exm-A-simple}
The algebra $A_2(r;\zeta)$ is simple if and only if $r=N$ and $\zeta\neq0$.
\end{proposition}

\begin{proof}
If $\zeta=0$, the ideal generated by $Y$ is a nonzero proper ideal of $A_2(r;0)$, so $A_2(r;0)$ is not simple.
If $\zeta\neq0$ and $r<N$, then $Y^r$ is a nontrivial central element, so $A_2(r;\zeta)$ is not simple.
Suppose that $r=N$ and $\zeta\neq0$.
Then $G$ and $Y$ satisfy $G^N=1$, $Y^N=\zeta$ and $GY=q^{-2}YG$ in $A_2(N;\zeta)$.
Since $q^{-2}$ is a primitive $N$-th root of unity, $A_2(N;\zeta)$ is isomorphic to $M_N(\bfk)$ as an algebra and is therefore simple.
\end{proof}

\begin{corollary}\label{cor:exm-A2-factorizability}
The reflective algebra $R_H(A_2(r;\zeta))$ is factorizable if and only if $r=N$ and $\zeta\neq0$.
In this case, the structure of $R_H(A_2(N;\zeta))$ as a left $H$-comodule algebra is obtained by setting $r=N$ in \cref{prop:exm-reflective-presentation}.
In particular, the relations for $g,y$ are $g^N=1$, $y^N=\zeta$ and $gy=q^{-2}yg$.
\end{corollary}

\begin{proof}
This follows from \cref{thm:reflective-factorizability-criterion,prop:exm-A-simple}.
\end{proof}

We now use $A_2(N;0)$ to give the example mentioned in \cref{rem:theta-bijective-factorizable,rem:comparison-WY-reflective}.

\begin{example}\label{ex:exm-reflective-theta-bijective-not-H-simple}
Let $A=A_2(N;0)$.
The algebra $A$ is generated by $G,Y$ subject to $G^N=1$, $Y^N=0$ and $GY=q^{-2}YG$.
In this case, $Z(A)=\bfk$.
On the other hand, $A$ is not simple by \cref{prop:exm-A-simple}.
Hence $\theta_{R_H(A)}:H^*\to E(H,R_H(A))$ is bijective by \cref{thm:theta-reflective-bijective-center}, but $R_H(A)$ is not $H$-simple by \cref{prop:factorizable-H-RH-simple}.
Therefore, $(R_H(A),K_{\mathrm{ref}}(A))$ is not factorizable.
\end{example}


\begin{thebibliography}{EGNO15}

\bibitem[AV22]{AV22}
A.~Appel and B.~Vlaar.
\newblock Universal {$K$}-matrices for quantum {K}ac--{M}oody algebras.
\newblock {\em Represent. Theory}, 26:764--824, 2022.

\bibitem[AV25]{AV25}
A.~Appel and B.~Vlaar.
\newblock Tensor {$K$}-matrices for quantum symmetric pairs.
\newblock {\em Comm. Math. Phys.}, 406:100, 2025.

\bibitem[BK19]{BK19}
M.~Balagovi\'{c} and S.~Kolb.
\newblock Universal {$K$}-matrix for quantum symmetric pairs.
\newblock {\em J. Reine Angew. Math.}, 747:299--353, 2019.

\bibitem[BZBJ18]{BZBJ18}
D.~Ben-Zvi, A.~Brochier, and D.~Jordan.
\newblock Quantum character varieties and braided module categories.
\newblock {\em Selecta Math. (N.S.)}, 24(5):4711--4748, 2018.

\bibitem[Bro13]{Bro13}
A.~Brochier.
\newblock Cyclotomic associators and finite type invariants for tangles in the solid torus.
\newblock {\em Algebr. Geom. Topol.}, 13(6):3365--3409, 2013.

\bibitem[Bur20]{Bur20}
S.~Burciu.
\newblock Representations and conjugacy classes of semisimple quasitriangular {H}opf algebras.
\newblock {\em SIGMA}, 16:039, 2020.

\bibitem[CSZ25]{CSZ25}
K.~Coulembier, M.~Stroi\'nski, and T.~Zorman.
\newblock Simple algebras and exact module categories.
\newblock {\tt arXiv:2501.06629}, 2025.

\bibitem[DN21]{DN21}
A.~Davydov and D.~Nikshych.
\newblock Braided {P}icard groups and graded extensions of braided tensor categories.
\newblock {\em Selecta Math. (N.S.)}, 27:65, 2021.

\bibitem[DKM03]{DKM03}
J.~Donin, P.~P.~Kulish, and A.~I.~Mudrov.
\newblock On a universal solution to the reflection equation.
\newblock {\em Lett. Math. Phys.}, 63(3):179--194, 2003.

\bibitem[EGNO15]{EGNO15}
P.~Etingof, S.~Gelaki, D.~Nikshych, and V.~Ostrik.
\newblock {\em Tensor categories}, volume~205 of {\em Mathematical Surveys
  and Monographs}.
\newblock American Mathematical Society, Providence, RI, 2015.

\bibitem[FRS02]{FRS02}
J.~Fuchs, I.~Runkel, and C.~Schweigert.
\newblock {TFT} construction of {RCFT} correlators {I}: Partition functions.
\newblock {\em Nuclear Phys. B}, 646:353--497, 2002.

\bibitem[Kas95]{Kas95}
C.~Kassel.
\newblock {\em Quantum groups}, volume~155 of {\em Graduate Texts in
  Mathematics}.
\newblock Springer-Verlag, New York, 1995.

\bibitem[Kol20]{Kol20}
S.~Kolb.
\newblock Braided module categories via quantum symmetric pairs.
\newblock {\em Proc. Lond. Math. Soc. (3)}, 121(1):1--31, 2020.

\bibitem[LWY25]{LWY25}
R.~Laugwitz, C.~Walton, and M.~Yakimov.
\newblock Reflective centers of module categories and quantum {$K$}-matrices.
\newblock {\em Forum Math. Sigma}, 13:e95, 2025.

\bibitem[LO17]{LO17}
S.~Lentner and T.~Ohrmann.
\newblock Factorizable {$R$}-matrices for small quantum groups.
\newblock {\em SIGMA}, 13:076, 2017.

\bibitem[Maj95]{Maj95}
S.~Majid.
\newblock {\em Foundations of quantum group theory}.
\newblock Cambridge University Press, Cambridge, 1995.

\bibitem[Mon93]{Mon93}
S.~Montgomery.
\newblock {\em Hopf algebras and their actions on rings}, volume~82 of {\em
  CBMS Regional Conference Series in Mathematics}.
\newblock American Mathematical Society, Providence, RI, 1993.

\bibitem[MW26]{MW26}
M.~M\"uller and C.~Walton.
\newblock Morita invariants of quasitriangular coideal subalgebras.
\newblock {\em J. Pure Appl. Algebra}, 230(3):108210, 2026.

\bibitem[NSS25]{NSS25}
D.~Nakamura, T.~Shibata, and K.~Shimizu.
\newblock Exact module categories over {$\mathrm{Rep}(u_q(\mathfrak{sl}_2))$}.
\newblock {\tt arXiv:2503.21265}, 2025.

\bibitem[Rad12]{Rad12}
D.~E.~Radford.
\newblock {\em Hopf algebras}, volume~49 of {\em Series on Knots and Everything}.
\newblock World Scientific Publishing Co. Pte. Ltd., Hackensack, NJ, 2012.

\bibitem[RT91]{RT91}
N.~Reshetikhin and V.~G.~Turaev.
\newblock Invariants of 3-manifolds via link polynomials and quantum groups.
\newblock {\em Invent. Math.}, 103(3):547--597, 1991.

\bibitem[Saf21]{Saf21}
P.~Safronov.
\newblock A categorical approach to quantum moment maps.
\newblock {\em Theory Appl. Categ.}, 37:818--862, 2021.

\bibitem[Shi19]{Shi19}
K.~Shimizu.
\newblock Non-degeneracy conditions for braided finite tensor categories.
\newblock {\em Adv. Math.}, 355:106778, 2019.

\bibitem[Skl88]{Skl88}
E.~K.~Sklyanin.
\newblock Boundary conditions for integrable quantum systems.
\newblock {\em J. Phys. A}, 21(10):2375--2389, 1988.

\bibitem[Skr07]{Skr07}
S.~Skryabin.
\newblock Projectivity and freeness over comodule algebras.
\newblock {\em Trans. Amer. Math. Soc.}, 359(6):2597--2623, 2007.

\bibitem[Tak79]{Tak79}
M.~Takeuchi.
\newblock Relative {H}opf modules---equivalences and freeness criteria.
\newblock {\em J. Algebra}, 60(2):452--471, 1979.

\bibitem[Tak94]{Tak94}
M.~Takeuchi.
\newblock Quotient spaces for {H}opf algebras.
\newblock {\em Comm. Algebra}, 22(7):2503--2523, 1994.

\bibitem[WY25]{WY25}
C.~Walton and H.~Yadav.
\newblock Nondegenerate module categories.
\newblock {\em Math. Z.}, 310:37, 2025.

\bibitem[WZ25]{WZ25}
J.~Wu and K.~Zhou.
\newblock Splitting property of quasitriangular {H}opf algebras.
\newblock {\tt arXiv:2501.15401}, 2025.
\end{thebibliography}
\end{document}